\documentclass[12pt]{amsart}
\usepackage{amsmath,amssymb,amsthm}
\usepackage[all]{xy}
\usepackage{xcolor}
\usepackage[utf8x]{inputenc}
\usepackage[T1]{fontenc}

\usepackage{amsthm}
\usepackage{listings}
\usepackage{graphicx}
\usepackage{hyperref}
\usepackage{cleveref}
\usepackage{caption}
\theoremstyle{theorem}
\newtheorem{theorem}{Theorem}[section]
\newtheorem{lemma}[theorem]{Lemma}
\newtheorem{cor}[theorem]{Corollary}
\theoremstyle{definition}
\newtheorem{definition}[theorem]{Definition}
\theoremstyle{remark}
\newtheorem{remark}[theorem]{Remark}

\numberwithin{equation}{section}
\begin{document}
\newcommand{\Exp}{\operatorname{Exp}}
\newcommand{\Mod}{\operatorname{Mod}}
\newcommand{\Inf}{\operatorname{Inf}}
\newcommand{\Hom}{\operatorname{Hom}}
\newcommand\PGL{\operatorname{PGL}}
\newcommand\GL{\operatorname{GL}}
\newcommand\tr{\operatorname{tr}}
\newcommand\lcm{\operatorname{lcm}}
\newcommand\FSexp{\operatorname{FSexp}}
\newcommand\Fun{\operatorname{Fun}}
\newcommand\Id{\operatorname{Id}}
\newcommand\CC{\mathbb C}
\newcommand\CCu{\CC^\times}
\newcommand\ZZ{\mathbb Z}
\newcommand\C{\mathcal C}
\newcommand\M{\mathcal M}
\newcommand\semdir{\rtimes}
\title[Pointed fusion categories of small rank]{Morita equivalence of pointed fusion categories of small rank}
\author{Micha\"el Mignard, Peter Schauenburg}

\begin{abstract}
We classify pointed fusion categories $\mathcal{C}(G,\omega)$ up to Morita equivalence for $1<|G|<32$. Among them, the cases $|G|=2^3 , 2^4$ and $3^3$ are emphasized. Although the equivalence classes of such categories are not distinguished by their Frobenius-Schur indicators, their categorical Morita equivalence classes are distinguished by the set of the indicators and ribbon twists of their Drinfeld centers. In particular, the modular data are a complete invariant for the modular categories $\mathcal{Z}(\mathcal{C}(G,\omega))$ for $|G[<32$. We use the computer algebra package GAP and present codes for treating complex-valued group cohomology and calculating Frobenius-Schur indicators.
\end{abstract}
\maketitle
\tableofcontents
\section{Introduction}
\label{sec:intro}
Fusion categories are semisimple rigid monoidal categories with finitely many isomorphism classes of simple objects, such that the space of automorphisms of every simple object is of dimension one. Classical examples are the category of representations of a finite group, and more generally the category of representations of finite dimensional semi-simple Hopf algebras; and hence the representations of quantum groups. Classical references on fusion categories include \cite{EGNO} and \cite{ENO}, and a more genereal review can be found in \cite{Mue}.

A first toy example of fusion categories is that of pointed fusion categories: A fusion category is called pointed if all simple objects are invertible. Therefore the set of simple objects has the structure of a finite group $G$, and the pentagon axiom for associativity leads naturally to a $3$-cocycle on the bar resolution of $G$. 

The dual of a fusion category $\C$ with respect to an indecomposable $\C$-module category $\M$ is the functor category $\C_{\M}^{\ast}:= \Fun_{\C}(\M,\M)$. A category $\mathcal{D}$ that is equivalent to the category $\C_\M^\ast$ is called categorically (or weakly) Morita equivalent to $\C$, \cite{Ni}. Categorical Morita equivalence of categories is an equivalence relation, see \cite{Mue2}. Moreover, two fusion categories are categorically Morita equivalent if and only if their Drinfeld centers are equivalent as braided fusion categories. The structure of a dual category over a pointed fusion category $\C(G,\omega)$ is given by the additional cohomological data $(H,\mu)$ where $H$ is a subgroup of $G$ and $\mu$ a $2$-cochain on $H$ such that $\partial \mu=\omega |_{H \times H \times H}$. The quadruple $(G,H,\omega,\mu)$ is called group-theoretical data and the associated dual category a group-theoretical category. Group-theoretical categories form a huge class of fusion categories, with interactions in many areas of both mathematics and physics.

One of these interactions arose from a reconstruction theory of orbifold vertex operator algebras studied in \cite{DPR}, in the form of twisted doubles of finite groups. A twisted quantum double $D_\omega(G)$ of a finite group $G$ is a certain quasi-triangular quasi-Hopf algebra defined by a $3$-cocycle on $G$. The study of these algebras and their representation categories, which are modular categories and in particular braided fusion categories, is important in several mathematical physics contexts. It is, in particular, interesting to classify such algebras up to so-called gauge equivalence. A definition of gauge transformation for quasi-bialebras can be found in \cite{Kassel}, but we will take a categorical point of view of gauge equivalence. Indeed, two doubles $D_\omega(G)$ and $D_{\omega'}(G')$ are gauge equivalent if and only if their module categories are equivalent as braided fusion categories. Also, the category $D_\omega(G)-\Mod$ is equivalent to the Drinfeld center $\mathcal{Z}(\C(G,\omega))$. Therefore, gauge equivalence classes of doubles $D_\omega(G)$ coincide with Morita equivalence classes of pointed fusion categories $\C(G,\omega)$. 

We note that \cite{GMN} and \cite{Uribe} completely classify Morita equivalent pointed categories for groups of order $8$, though both with different points of view and techniques. We generalize this result to all groups whose order lies between $2$ and $31$. We mainly use two ingredients in order to obtain such a classification. The first one is the main tool of \cite{GMN} that states a particular Morita equivalence between $\C(G,\omega)$ and $\C(G',\omega')$ for special pairs $(G,\omega)$ and $G',\omega)$. To wit: $G$ and $G'$ are extensions $1 \to A \to G \to Q \to 1$ and $1 \to \hat{A} \to G' \to Q \to 1$ with $A$ abelian, and the cocycles $\omega, \omega'$ can be computed by inflations of cocycles on the quotient and some extra data. The second one is an explicit computation of the Frobenius-Schur indicators of pointed fusion categories, their dual categories and their centers; we use also the ribbon twists of these centers.

Frobenius-Schur indicators of pivotal categories are generalisations of the second indicator $\nu_2(\chi)$ defined for an irreducible character of a finite group by its namesakes in $1906$. Several generalisations of those indicators, including categorical second indicators \cite{Bantay}\cite{Fuchs}, higher indicators for modules over semi-simple Hopf algebras \cite{LinMont} and second indicators for modules over semi-simple quasi Hopf algebras \cite{MN2}, led to the definition of higher indicators in \cite{NS1}. A new formula for higher indicators of general group-theoretical categories has been found in \cite{S1}, including formulas for pointed categories and modules over a twisted quantum double of a finite group. The proof of this result is strongly based on a trick used by \cite{Na} that permits us to change a general modular data to an equivalent one which can be described with restricted cohomology settings. These formulas will allow us to compute the indicators with GAP based programms.

In fact we view it as a goal of the paper of independent interest to implement the computation of indicators in GAP. Also, the construction of group-theoretical categories requires making the cohomological data involved in this construction available in an efficient fashion, and we think that the approach to this problem we will present is also of independent interest.

The paper is organized as follows: Section \ref{sec:GAPcoho} gathers preliminary results concerning calculations with group cohomology wich is useful to describe both the objects we are dealing with and their Frobenius-Schur indicators. In Section \ref{sec:context}, general background about pointed fusion categories is recalled and we give formulas for the Frobenius-Schur indicators of these categories. Our results are given in section \ref{sec:work} and can be summarized by \cref{completeclass,completeinv}. 

Throughout the paper, linear categories will be taken over the field of complex numbers, and all groups are assumed to be finite.
\section{Complex valued cohomology with GAP}
\label{sec:GAPcoho}
In this section, we present GAP/HAP code that calculates the cohomology groups of a finite group $G$ with coefficients in $\mathbb C^\times$. We take for granted that we have suitable resolutions available (usually provided by HAP).

It is well-known (see below) that $H^n(G,\mathbb C^\times)$ is a finite abelian group, but it is defined as the cohomology group of a cochain complex of abelian groups that are not finitely generated (simply since $\mathbb C^\times$ is already not finitely generated as an abelian group). Thus, even if the cohomology group we need is finite, we cannot naively compute it by feeding the complex defining it to a computer algebra system.

Denote by $\mu_m\subset \mathbb C^\times$ the subgroup of $m$-th roots of unity. As we will discuss below in detail, the Universal Coefficient Theorem implies that $H^n(G,\mathbb C^\times)$ is naturally an epimorphic image of $H^n(G,\mu_m)$ for suitable $m$. The latter group now is defined as the cohomology group of a cochain complex of \emph{finite} abelian groups, and thus readily accessible to machine calculations. In other words, every cohomology class in $H^n(G,\mathbb C^\times)$ is represented by a cocycle taking values in $\mu_m$. This is the way proposed by the makers of HAP to obtain all cohomology classes in $H^n(G,\mathbb C^\times)$ \cite{hap}.

A slight drawback is that the kernel of the map $H^n(G,\mu_m)\to H^n(G,\mathbb C^\times)$ is not trivial. In other words, a nontrivial cocycle with values in $\mu_m$ may be trivial when regarded as a cocycle with values in $\mathbb C^\times$. While the kernel is known through the Universal Coefficient Theorem and not hard to compute as an abstract group (it is the Ext-group of a cyclic group by a homology group), it is not immediately available as a subgroup of $H^n(G,\mu_m)$.

We use the following approach. The Universal Coefficient Theorem implies that $H^n(G,\mathbb C^\times)$ is isomorphic to $\operatorname{Hom}(H_n(G),\mu_m)$ for suitable $m$, and $H_n(G)$ can be obtained as the homology of a complex of finitely generated free abelian groups. The isomorphism and the form of all homomorphisms defined on $H_n(G)$ can be made sufficiently explicit as to obtain an (additive) presentation
\[H^n(G,\mathbb C^\times)=\langle x_1,\dots,x_k|s_ix_i=0\rangle,\]
i.~e.\ an explicit isomorphism
\[H^n(G,\mathbb C^\times)\cong \mathbb Z/(s_1)\times\dots\times\mathbb Z/(s_k).\]

We start with standard definitions. Let K be a chain complex of abelian groups $K_n$, with differential maps $\partial_n : K_n \rightarrow K_{n-1}$.
Denote by $Z_n(K),B_n(K)\subset K_n$ the subgroups of cycles, resp.\ boundaries, 
and by $H_n(K):=Z_n(K)/B_n(K)$ the $n$-th homology group.

For an abelian group $A$ we consider the cochain complex
$Hom(K,A)$ with differential 
$\delta^n (f) := f \circ \partial_{n+1}$. We denote by 
$Z^n(K,A),
B^n(K,A),$
and $H^n(K,A)$
the groups of cocycles, coboundaries, and the cohomology group of this cochain complex.

We recall the following standard result :
\begin{theorem}\label{UCTtheo}
(Universal Coefficient Theorem)

Let $K$ be a complex of free abelian groups. 
One has a short exact sequence of abelian groups:
\begin{equation*}
E_A \ : \  0 \rightarrow Ext(H_{n-1}(K),A) \overset{\beta_A}\hookrightarrow H^n(K,A) \overset{\alpha_A} \twoheadrightarrow Hom_{\mathbb{Z}}(H_n(K),A) \rightarrow 0
\end{equation*}
This sequence is natural in both K and A. Moreover, it has a (non-natural) splitting 
$\pi_A\colon Hom_{\mathbb{Z}}(H_n(K),A) \rightarrow H^n(K,A)$
\end{theorem}
We use this theorem in the following way. For any positive integer $m$, define the map $\phi : \mathbb{Z}/m \mathbb{Z} \rightarrow \mathbb{C}^*, n \mapsto \zeta_m ^n$, where $\zeta_m$ denotes a primitive $m$-th root of unity. Since the exact sequence given by the UCT is natural in $A$, we have a chain map between the exact sequences $E_{\mathbb{C}^*}$ and $ E_{\mathbb{Z}/m \mathbb{Z}}$. Also, since $\mathbb{C}^*$ is divisible, $Ext(C,\mathbb{C}^*)$ is trivial for any abelian group $C$. We obtain then the following diagram with exact rows:
{\tiny
\begin{displaymath}
\xymatrix{
   & 0 \ar[r]  & H^n(K,\mathbb{C}^\times) \ar[r]^-{\cong} & Hom_{\mathbb{Z}}(H_n(K),\mathbb{C}^\times) \ar[r]  &0\\
    0 \ar[r] & Ext(H_{n-1}(K),\mathbb{Z}/m \mathbb{Z}) \ar[u] \ar[r] & H^n(K,\mathbb{Z}/m \mathbb{Z}) \ar[u]^-{\gamma} \ar@{->>}[r]  & Hom_{\mathbb{Z}}(H_n(K),\mathbb{Z}/m \mathbb{Z}) \ar[u] \ar[r] & 0
    }
\end{displaymath}}
Put $m:=\Exp(H_n(K))$. Then the rightmost vertical arrow from $Hom_{\mathbb{Z}}(H_n(K),\mathbb{C}^\times)$ to $Hom_{\mathbb{Z}}(H_n(K),\mathbb{Z}/m \mathbb{Z})$ is an isomorphism. By naturality, the right hand square commutes, and so the map $\gamma$ from $H(K,\mathbb{Z}/m \mathbb{Z})$ to $H^n(K,\mathbb{C}^\times)$ is surjective. Moreover, the composition $\gamma \circ \pi_{\mathbb{Z}/m \mathbb{Z}}$ is an isomorphism from $Hom_{\mathbb{Z}}(H_n(K),\mathbb{Z}/m \mathbb{Z})$ to $H^n(K,\mathbb{C}^\times)$.

We now show how to implement this in GAP. All codes are gathered in \ref{ssec:GAP1}.
\subsection{Calculating cohomology groups}
\label{ssec:=UCT}
We describe in detail how we compute cohomology groups using GAP/HAP code. From now on $K$ will always be a complex of \emph{free} $\mathbb{Z}$-modules. Also, we  fix a basis ${\mathcal{K}}:=(e_1,\dots,e_N)$ for $K_n$, where $N$ is its rank. Cochains, coboundaries and boundary maps are then  described by integer row vectors and matrices. The strategy is to use the isomorphism between $Hom_{\mathbb{Z}}(H_n(K),\mathbb{Z}/m \mathbb{Z})$ and $H^n(K,\mathbb{C}^\times)$ instead of just the surjection $\gamma$: we compute the group $Hom_{\mathbb{Z}}(H_n(K),\mathbb{Z}/m \mathbb{Z})$, and then lift its elements to $H^n(K,\mathbb{C}^\times)$, using the composition of a section of the UCT sequence for $\mathbb{Z}/m \mathbb{Z}$ and the surjection from $H^n (K,\mathbb{Z}/m \mathbb{Z})$ to $H^n(K,\mathbb{C}^\times)$ as described before.

First, we give details about the section map of the UCT. The notations are the same as in  \cref{UCTtheo}. Following \cite{Hatcher} ---where the reader will find a complete proof of \cref{UCTtheo}--- we construct the map $\alpha_A$ and its right inverse. 

Take $f \in Z^n(K,A)$. The cohomology class of $f$ in $H^n(K,A)$ is denoted by $\bar{f}$. By definition, we have that $\delta^n(f)=f \circ \partial_{n+1}=0$. So $f|_{B_n(K)}=0$ and $f$ factors over $H_n(K)$. Take $\alpha_A(\bar{f}):=f|_{Z_n(K)}$. We have to verify that $\alpha_A$ is actually well-defined. If $f \in B^n(K,A)$, then there exists $\tilde{f}: K_{n-1} \rightarrow A$ such that $f=\delta^{n-1}(\tilde{g})=\tilde{f} \circ \partial_n$. Therefore, $f \circ \partial_{n+1}=0$ and $\alpha_A$ is well-defined.

Now, the exact sequence $0 \rightarrow Z_n(K) \hookrightarrow K_n  \twoheadrightarrow B_{n-1}(K) \rightarrow 0$ splits since $B_{n-1}(K)$ is free. Let $\pi$ denote a splitting of the injection $Z_n(K) \hookrightarrow K_n$, that is $\pi$ is a surjective map from $K_n$ to $Z_n(K)$ which restricts to the identity on $Z_n(K)$. Composing $\pi$ with a function $f : Z_n(K) \rightarrow A$ that vanishes on $B_n$ gives a function $\tilde{f}:=f \circ \pi : K_n \rightarrow A$ that still vanishes on $B_n$. Also, we have that $\alpha (\tilde{f}) = \alpha(f \circ \pi)=(f \circ \pi)|_{Z_n(K)}=f \circ id_{Z_n(K)}=f$, so the composition with $\pi$ gives a splitting of $\alpha$, denoted $\pi_A$.

Now we get some information about $H_n (K)$, which is the quotient of the free modules $Z_n (K)$ and $B_n (K)$. Let $k,l$ be their ranks, respectively, and we assume that we also have bases ${\mathcal{Z}}:=(z_1 ,\dots , z_k)$ and ${\mathcal{B}}:=(b_1 , \dots ,b_l)$. We recall the following theorem for free modules over euclidean rings :
\begin{theorem}[Adapted Basis Theorem]\label{abt}
Let $A$ be an euclidean ring. If $F$ is a free module over $A$ of rank $k$, and $S$ is a free submodule of $F$ of rank $l$, then there exists a basis ${\mathcal{B}}=\{ b_1 , \dots , b_k\}$ of $F$ and non-zero elements  $s_1, \dots , s_l$ of $A$, such that $s_1 | s_2 | \dots | s_l$ and the set $\{ s_1 \cdot b_1, \dots ,s_l \cdot b_l \}$ is a basis of $S$.
\end{theorem}
This theorem is equivalent with the existence of the Smith Normal Form for matrices with coefficients in $A$. To wit: if one has a basis for $F$ and a system of generators for $S$, take $M$ the matrix of decomposition of this system of generators into the basis of $F$. Then the elements $s_i$, called the invariant factors, are the coefficients of the SNF of M. There exists in GAP a function that computes such bases from provided systems of generators for both $F$ and $S$. We slightly complete it, so that it also gives the change of bases needed for our implementation.

Remember that $Z_n(K)$ and $B_n(K)$ are submodules of $K_n$ with basis $(e_1,\dots , e_N)$, which we identify as the free module $\mathbb{Z}^N$ with its canonical basis.  For a system $\Gamma=(\gamma_1 ,\dots \gamma_k)$ of elements in $K_n$, we will also denote as $\Gamma$ the matrix formed by the decomposition as row vectors of the $\gamma_i$ in the canonical basis.

We obtain then, from a basis ${\mathcal{Z}}:=\{z_1 ,\dots , z_k\}$ of $Z_n(K)$, and a basis (or system of generators) ${\mathcal{B}}$ of $B_n(K)$, a new basis ${\mathcal{Z}}':=\{z'_1 ,\dots , z'_k\}$ of $Z_n(K)$ and the invariant factors $\{s_1 , \dots ,s_l\}$ such that $\{s_1 \cdot b_1 ,\dots ,s_l \cdot b_l \}$ is a basis of $B_n(K)$. We denote $C$ the matrix such that ${\mathcal{Z}}'=C {\mathcal{Z}}$.

Therefore, the quotient $H_n(K)=Z_n(K) / B_n(K)$ is a finite abelian group isomorphic to $\mathbb{Z}/s_1 \mathbb{Z} \times \cdots \times \mathbb{Z}/s_l \mathbb{Z}$. Also, as $s_1 | \dots | s_l$ we have that the exponent of this group is $s_l$.
Now, we want a description of the set $Hom_{\mathbb{Z}}(H_n(K),\mathbb{Z}/m \mathbb{Z})$ where $m$ is (a multiple of) the exponent $s_k$. The elements of this set are homomorphisms $f$ from $Z_n (K)$ to $\mathbb{Z}/m \mathbb{Z}$ such that $f|_{B_n (K)}=0$. For a morphism $f : Z_n (K) \rightarrow \mathbb{Z}/m \mathbb{Z}$, we denote $\bar{f}:=(f(z_1'), \dots , f(z_k'))$. This row vector represents then the morphism $f$ in the basis ${\mathcal{Z}}'$.

For all $ 1 \leqslant i \leqslant l$, we have :
$$f(b_i')=f(s_i \cdot z_i ')=s_i \cdot f(z_i ')$$
Then, $f|_{B_n (K)}=0 \Leftrightarrow s_i \cdot f(z_i ')=0$ mod $m$ for all $1 \leqslant i \leqslant l$.

One can then describe an explicit system of generators for 
\\$\Hom_{\mathbb{Z}}(H_n(K),\mathbb{Z}/m \mathbb{Z})$ by setting, for $1 \leqslant i \leqslant l$, the vector $\bar{h_i}$ of size $k$: 
$$\bar{h_i}:=\left(0,\dots,\frac{m}{s_i},\dots,0\right)$$

where the only non-zero coordinate is in the $i$-th position.

The morphisms $h_i$  represented in ${\mathcal{Z}}'$ by the row vectors $\bar{h_i}$ are then a system of generators of $Hom_{\mathbb{Z}}(H_n(K),\mathbb{Z}/m \mathbb{Z})$. To express them in the basis ${\mathcal{Z}}$, one has just to multiply those vectors $\bar{h_i}$ by the matrix $C^{-1}$ (and reduce the results modulo $m$) ; we denote those new row vectors by $\tilde{h_i}$.
Just note the fact that if $s_i$ is equal to $1$, the corresponding vector $\tilde{h_i}$ is equal to the null vector modulo $m$. Once those removed, each remaining vector represents a cyclic part of the abelian group $Hom_{\mathbb{Z}}(H_n (K), \mathbb{Z}/m \mathbb{Z})$. Thus, the list of orders of those vectors modulo $m$, which is the list of the $s_i$ that are not equal to $1$, gives us the torsion of $H^n(K,\mathbb{C}^\times)$.

Then we lift those generators $h_i$ to elements of $H^n (K,\mathbb{Z}/m \mathbb{Z})$. For this, we are going to use the section of the UCT we described earlier.
 
We use \cref{abt} for $K_n$ and $Z_n(K)$ ; the obtained basis of $K_n$ is denoted ${{K}}'=\{\kappa_1 , \dots , \kappa_N\}$. The advantage here is that, as the quotient of $K_n$ by $Z_n(K)$ is isomorphic to the free module $B_{n-1}(K)$, the invariant factors are all equal to $1$, and so ${\mathcal{Z}}'':=\{\kappa_1 ,\dots ,\kappa_k\}$ is a basis for $Z_n(K)$. We denote respectively $C_1$ and $R_1$ the matrices such that ${\mathcal{Z}}''=C_1 {\mathcal{Z}}$ and ${\mathcal{K}}'=R_1 {\mathcal{K}}$.

Thus, one can easily express in those bases an $(N \times k)$-matrix $I$ coding a surjective map from $K_n$ to $Z_n(K)$, denoted $\pi$, such that :
$$ \pi( \kappa_i) = \kappa_i\ ,\ 1 \leqslant i \leqslant k$$ 
$$ \pi( \kappa_i)=0 \ ,\ k+1 \leqslant i \leqslant N$$
Then, for a vector $\tilde{h}$ representing a morphism $h:H_n(K) \rightarrow \mathbb{Z}/m\mathbb{Z}$ in the basis ${\mathcal{Z}}$, the vector $\hat{h}:=C_1 \times \tilde{h}$ represents $h$ in the basis ${\mathcal{Z}}''$. Applying the matrix $I$ to $\hat{h}$ gives a row vector, denoted $\widehat{\omega}$, with $N$ coordinates, whose are the image in $\mathbb{Z}/m \mathbb{Z}$ of the basis elements $\{\kappa_1 ,\dots ,\kappa_N \}$ by $h \circ \pi$. The application $\omega : K_n \rightarrow \mathbb{Z}/m \mathbb{Z}$ is an $n$-cocycle for the complex $K$. Finally, multiplying this row vector $\widehat{\omega}$ by the matrix $R_1^{-1}$ gives a row vector $\tilde{\omega}$ that represents $\omega$ in the canonical basis of $K_n$.

We obtain by this procedure, from the set of vectors $\tilde{h_i}$, a set of vectors $\tilde{\omega_i}$ representing elements $\omega_i=h_i \circ \pi$ of $H^n (K,\mathbb{Z}/m \mathbb{Z})$. Using the fact that the composition of the section of the UCT and the (surjective) map $\gamma$ is an isomorphism, each $\gamma (\omega_i)$ generates a different cyclic part of order $s_i$ in $H^n (K,\mathbb{Z}/m \mathbb{Z}) \cong \mathbb{Z}/s_1 \mathbb{Z} \times \cdots \times \mathbb{Z}/s_l \mathbb{Z}$.
\subsection{Functoriality}
\label{ssec:functor}
Now, we consider two complexes $K$ and $K'$ and a chain map $f\colon K\to K'$. We keep the same notations as in the previous section, and add a single quote to all objects associated to $K'$. We have seen how to compute $H^n(K,\CC^\times)$ explicitly through an isomorphism with a group of homomorphisms $\Hom(H_n(K),\ZZ/m\ZZ)$. Now we want to extend this description to cover the induced homomorphism $H^n(f,\CC^\times)$ by completing the diagram
  \begin{equation*}
    \xymatrix{H^n(K',\CCu)\ar[rr]^{H^n(f,\CCu)}\ar[d]^{\cong}&&H^n(K,\CCu)\ar[d]^{\cong}\\
             \Hom(H_n(K'),\ZZ/m'\ZZ))\ar[rr]^{{\hat{F}}}&&\Hom(H_n(K),\ZZ/m\ZZ))}
  \end{equation*}

We get first a matrix description of the homomorphim $H_n(f): H_n(K) \to  H_n(K')$. From the chain map $f : K \rightarrow K'$ we form the $(N \times N')$-matrix $F$ whose rows are the images of the canonical basis of $K_n$ by $f_n$, decomposed in the canonical basis of $K'_n$. Remember that a morphism $h':H_n(K') \to \mathbb{Z}/m' \mathbb{Z}$, where $m'$ is a multiple of $\Exp(H_n(K'))$, is described by a row vector of size $N'$ with integral coefficients, but whose supposed to been seen modulo $m'$.

However, as $m'$ was taken to be a multiple of the exponent of $H_n(K')$ and because in general this is not a multiple of the exponent of $H_n(K)$, we cannot just multplying this vector by $F$ in order to get the composition $H_n(f) \circ h'$. We need first to extend the rings $\mathbb{Z}/m\mathbb{Z}$ and $\mathbb{Z}/m'\mathbb{Z}$ in such a way that this product makes sense. This procedure is summarized by the following diagram.

Take $p$ a multiple of $\lcm(m,m')$, $\xi_m$, $\xi_{m'}$ and $\xi_p$ respectively $m$, $m'$ and $p$ primitive roots of unity, $t$, $t'$ and $s$ generators of (respectively) $\mathbb{Z}/m\mathbb{Z}$, $\mathbb{Z}/m'\mathbb{Z}$ and $\mathbb{Z}/p\mathbb{Z}$, we construct $h : H_n(K) \rightarrow \mathbb{Z}/m\mathbb{Z}$ such that $\tilde{h'} \circ H_n(f) = \tilde{h}$, where $\tilde{h'}=(t' \mapsto \xi_{m'}) \circ h'$ and $\tilde{h}=(t \mapsto \xi_m) \circ h$.
\begin{displaymath}
\xymatrix{H_n(K') \ar@*{}[r]^{h'} \ar[rd]^{\tilde{h'}} & \mathbb{Z}/m'\mathbb{Z} \ar[d]^{(t' \mapsto \xi_{m'})} \ar@*{}@{^{(}->}[drr] \\
 & \mathbb{C}^\times &  &  \mathbb{Z}/p \mathbb{Z} \ar[ll]_{(s' \mapsto \xi_p)}\\
 H_ n(K) \ar@*{}[uu]^{H_n(f)} \ar@{-->}[r]^{h} \ar@{-->}[ur]^{\tilde{h}} & \mathbb{Z}/m \mathbb{Z} \ar[u]_{(t \mapsto \xi_m)} \ar@{^{(}->}[urr]
}
\end{displaymath}
Therefore, we multiply the coordinates of $h$ by the integer $\frac{p}{m}$, apply to the result the matrix $F$, and finally divide the coordinates of the formed row vector (of size $N$) by the integer $\frac{p}{m'}$. That the last vector representing $h$ is integral is garanted by the commutativity of the diagram. Reducing this vector modulo $m$ gives the image of $h'$ by composition with $H_ n(f)$, that is $\hat{F}(h')=h$.

Finally, in order to get a description of $H^n(f,\mathbb{C}^\times)$, we recall the isomorphism $\gamma \circ \pi_{\mathbb{Z}/m\mathbb{Z}} : \Hom(H_n(K),\mathbb{Z}/m\mathbb{Z}) \to H^n(K,\mathbb{C}^\times)$. That is, from a cocycle $\widetilde{\omega}' \in Z^n(K',\mathbb{Z}/m'\mathbb{Z})$ representing $\omega' \in Z^n(K',\mathbb{C}^\times)$, we form the homomorphism $h' \in \Hom(H_n(K'),\mathbb{Z}/m'\mathbb{Z})$ image of $\widetilde{\omega}'$ by the surjection of the UCT for $K'$. We compute then a homomorphism $h \in \Hom(H_n(K),\mathbb{Z}/m\mathbb{Z})$ be the previous procedure. Applying $\gamma \circ \pi_{\mathbb{Z}/m\mathbb{Z}}$ to $h$ gives a cocycle $\omega \in H^n(K,\mathbb{C}^\times)$ such that its cohomology class is the image of the cohomology class of $\omega'$ by $H^n(f,\mathbb{C}^\times)$. 
\begin{remark}
  In our applications, the complexes $K,K'$ will arise from tensoring free resolutions for groups $G,G'$ with $\ZZ$, and the chain map $f$ will arise from a group homomorphism $\phi\colon G\to G'$. The fact that we can represent cohomology classes in $H^n(K,\CCu)=H^n(G,\CCu)$ by homomorphisms $H_n(K)\to\ZZ/m\ZZ$ (and not by equivalence classes of cocycles) has computational advantages: For example, it is easier to take orbits of a cohomology group under the automorphism group of the group in question, if we act on homomorphisms instead of on equivalence classes modulo coboundaries. Similarly, it is easy to check whether the restriction of a cohomology class to a subgroup is trivial, or if a cohomology class is the inflation of another cohomology class with respect to a group epimorphism.

There is a caveat, however: While we can also obtain in this fashion a cocycle $\omega\colon K\to\ZZ/m\ZZ$ representing the image in $H^n(G,\CCu)$ of the cohomology class in $H^n(G',\CCu)$ represented by a given cocycle $\omega'\colon K'\to \ZZ/m'\ZZ$, it is generally not true that $\omega$ would be the restriction of $\omega'$ along $H_n(f)$ (even if we properly treat the difference between $m$ and $m'$). This is due to the fact that we had to choose a section $\pi_{\ZZ/m\ZZ}$ above. In particular, if we find that the restriction of a cohomology class to a subgroup is trivial, this does not give us a representative of the cohomology class whose restriction would be the trivial cocycle; the restriction will in general only be a coboundary, and we have no information on the cochain whose coboundary it is. In particular, if we pass from cocycles given in an arbitrary free resolution to cocycles $\omega\colon G^n\to\CCu$ with respect to the bar resolution, we may know that such a cocycle represents a cohomology class whose restriction to a subgroup $H\subset G$ is trivial, but we will not have $\nu\colon H^{n-1}\to\CCu$ in hand such that the restriction of $\omega$ is $d\nu$. We believe that it is not a problem to add such a feature to our treatment of $\CCu$-valued cohomology; clearly it would be useful since, for example, treating group-theoretical categories in full generality one has to computationally deal with just such a situation. However, since for our concrete problems we happened to chance upon sufficiently many cocycles where the restrictions in question are not only cohomologically, but rather outright trivial, we postponed this.
\end{remark}

\subsection{GAP implementation}
\label{ssec:GAP1}
We give now the codes we have implemented in order to compute the $\mathbb{C}^\times$-valued cohomology groups of finite groups, and the morphism $H^n(\phi):H^n(H,\mathbb{C}^\times) \rightarrow H^n(G,\mathbb{C}^\times)$ for a group morphism $\phi:G \rightarrow H$.
Let $G$ be a finite group. Let $R$ be a resolution of $\mathbb{Z}[G]$-modules, and $A$ a $\mathbb{Z}$-module. By the Hom-Tensor adjunction, we have that 
\begin{eqnarray*}
Hom_{\mathbb{Z}[G]} (R,A) & \cong & Hom_{\mathbb{Z}[G]} (R,Hom_{\mathbb{Z}}(A,\mathbb{Z}))\\
& \cong & Hom_{\mathbb{Z}}(R \otimes \mathbb{Z},A)
\end{eqnarray*}
HAP provides free resolutions of a given size $n$ for finite groups, by using for example the command
\begin{lstlisting}
R:=ResolutionFiniteGroup(G,n);
\end{lstlisting}
or other commands that may be available for specific groups; we note that it is essential to have a command that produces a free resolution with a contracting homotopy which HAP uses to calculate cocycles in the bar resolution.

Also, we can construct the chain complex $K:=R \otimes \mathbb{Z}$ by :
\begin{lstlisting}
K:=TensorWithIntegers(R);
\end{lstlisting}
Now, we can implement the algorithm previously described to obtain the torsion of the homology group $H^n(K,A)$, together with a system of generators --- one for each cyclic part.
\begin{lstlisting}[linerange=cohomology-end]
LoadPackage("hapco");


InstallMethod( OneMutable,true,[IsCcElement],0,c -> One(InCcGroup(c)) );




ComplementIntMatWithTransforms:=function( full,sub )
local F, S, M, r, T, c, R;
  F := BaseIntMat( full );
  if IsEmpty( sub ) or IsZero( sub ) then
    return rec( complement := F, sub := [], moduli := [] );
  fi;
  S := BaseIntersectionIntMats( F, sub );
  if S <> BaseIntMat( sub ) then
    Error( "sub must be submodule of full" );
  fi;
  # find T such that BaseIntMat(T*F) = S
  M := Concatenation( F, S );
  T := NormalFormIntMat( M, 4 ).rowtrans^-1;
  T := T{[Length(F)+1..Length(T)]}{[1..Length(F)]};
  c := NormalFormIntMat(T*F,6)!.rowtrans;

  # r.rowtrans * T * F = r.normal * r.coltrans^-1 * F
  r := NormalFormIntMat( T, 13 );
  M := r.coltrans^-1 * F;
  R := rec( complement := BaseIntMat( M{[1+r.rank..Length(M)]} ),
            sub := r.rowtrans * T * F,
            moduli := List( [1..r.rank], i -> r.normal[i][i] ) ,
	    fulltrans := r.coltrans^(-1),
	    subtrans := r!.rowtrans*c^(-1));
  return R;
end;


#& cohomology #&

UniversalCoefficientsTheorem:=function(K,n)
local ZZ,     # a basis for Z_n(K)
      B,      # a basis for B_n(K)
      r,      # the ABT for Z_n(K) and B_n(K)
      C,      # the matrix C
      D,      # the list of the invariant factors
      m,      # the exponent of the H_n(K)
      homlist,# a "basis" for Hom(H_n(K),Z/mZ)
      r2, R1,C1,
              # the different variables implicated
              # in the ABT for K_n and Z_n(K)
		
      I,      # the matrix I
      cocyclelist;
              # a "basis" for H^n(K,C^*)


	
ZZ:=NullspaceIntMat(TransposedMat(
                BoundaryMatrix(K,n)));
if not Size(ZZ)>0 then return 
	rec(	complex:=K,
		cyclebasis:=ZZ,
		cycletrans:=[],
		hombasis := [],
		cocyclebasis := [],
		torsion := [1],
		exponent := 1,
		lift:=[]
	);
fi;
B:=BaseIntMat(TransposedMat(
                 BoundaryMatrix(K,n+1)));
r:=ComplementIntMatWithTransforms(ZZ,B);
C:=r!.fulltrans;
D:=r!.moduli;
m:=D[Size(D)];
homlist:=[];
List(Filtered(D,x->not x=1),x->
	ListWithIdenticalEntries(Size(ZZ),0));
for i in [1..Size(D)] do
	if D[i]=1 then continue;fi;
	hom:=ListWithIdenticalEntries(Size(ZZ),0);
	hom[i]:=m/D[i];
	Add(homlist,hom);
od;
homlist:=List(homlist,x->((C^(-1))*x) mod m);

r2:=ComplementIntMatWithTransforms(
              IdentityMat(K!.dimension(3)),ZZ);
C1:=r2!.subtrans;
R1:=r2!.fulltrans;

I:=Concatenation(IdentityMat(Size(ZZ)),
     NullMat(K!.dimension(n)-Size(ZZ),Size(ZZ)));
cocyclelist:=List(homlist,x->(R1^(-1)*(I*(C1*x))) 
                     mod m);
	
return( rec(
	complex:=K,
	cyclebasis:=ZZ,
	cycletrans:=C,
	hombasis := homlist,
	cocyclebasis := cocyclelist,
	torsion := Filtered(D,x->not x=1),
	exponent := m,
	lift:=R1^(-1)*I*C1
		)
	);
end;


#& end #&


#& order #&

CocycleOrder:=function(f,m)
	local list;
	list:=List(Filtered(f,x->not x=0),
                   x->m/Gcd(x,m));
	if list=[] 
		then return 1;
		else return Maximum(list);
	fi;
end;    

#& end #&   

#& functoriality #&

GroupCohomologyFunctor:=function(K1,K2,phi,n)
local 	UCT1, 	# UCT for K1
	m1, 	# the exponent of H_n(K1)
	Z1,	# the adapted basis of Z_n(K1)
	UCT2, 	# UCT for K2
	m2, 	# the exponent of H_n(K2)
	Z2,	# the adapted basis of Z_n(K2)
	p, 	# p=lcm(m1,m2)
	F, 	# the matrix F
	hphi; 	# the morphism H_n(phi)
		

UCT1:=UniversalCoefficientsTheorem(K1,n);;
m1:=UCT1!.exponent;;
Z1:=UCT1!.cyclebasis;
s:=NormalFormIntMat(Z1,13);
UCT2:=UniversalCoefficientsTheorem(K2,n);;
m2:=UCT2!.exponent;;
Z2:=UCT2!.cyclebasis;
F:=List(Z1,
	x->SolutionIntMat(Z2,phi!.mapping(x,n)));;
p:=Lcm(m1,m2);;

hphi:=function(f)
	return ((F*(m1*f))/m2) mod m1;
end;
return rec(
	matrix:=F,
	mapping:=hphi,
	p:=p,
	m1:=m1,
	m2:=m2);
end;

#& end #&


#& proj #&

ProjectiveCharacters:=function(G,alpha,o)
local A,r,Ga,Gaa,ire,t,lift,testpro;

# we represent the cyclic group as a G-Outer group

A:=GOuterGroup();
SetActingGroup(A,G);
if o=1 then
    SetActedGroup(A,GroupByGenerators(
	[One(CyclicGroup(o))]));
    else
        SetActedGroup(A,CyclicGroup(o));
fi;
SetOuterAction(A,function(g,x) return x; end);

# we represent the cocycle f as a 
# standard 2-cocycle in HAP 

r:=Standard2Cocycle();
SetActingGroup(r,A!.ActingGroup);
SetMapping(r,function(x,y) return 
	(A!.ActedGroup.1)^alpha(x,y);
	end);
SetModule(r,A);


# now, we construct the extension of G by A along r
# we also give its character table

Ga:=CcGroup(A,r);
Gaa:=UnderlyingGroup(Ga);
ire:=Irr(Gaa);

# we give the element in Gaa that corresponds to 
# the generators of the cyclic group, and a 
# section of the epi of Gaa on G


t:=CcElement(FamilyObj(One(Ga)),GeneratorsOfGroup(
   Ga!.Fibre)[1],One(Ga!.Base),InCcGroup(One(Ga)));
lift:=function(g)
      return CcElement(FamilyObj(One(Ga)),
	One(Ga!.Fibre),g,InCcGroup(One(Ga)));
end;

# finally, we give the function that tests for a 
# charater of Gaa if it comes from a projective 
# character of G

testpro:=function(l)
	local list,chi;
	list:=[];
	for chi in l do
      		if t^chi/E(o)=(t^Size(Ga))^chi 
			then Add(list,chi); 
		fi;
	od;
	return list;
end;

return rec(
	Gaa:=Gaa,
	lift:=lift,
	table:=testpro(ire)
	);
end;

#& end #&


#& automorphisms #&

AutomorphismOrbits:=function(G,R,K,UCT)
local 	rcv,
	coho,
	aut,
	gens,
	autmat,
	orbits;

rcv:=function(x,n) 
     return List(x,y->ZmodnZObj(y,n));
     end;
coho:=NearAdditiveGroup(
      List(UCT!.hombasis,x->rcv(x,UCT!.exponent)));
aut:=AutomorphismGroup(G);
gens:=GeneratorsOfGroup(aut);
autmat:=List(
         gens,phi->TransposedMat(
          GroupCohomologyFunctor(
           K,K,TensorWithIntegers(
	    EquivariantChainMap(R,R,phi)),
           3)!.matrix));
autmat:=List(
         autmat,x->List(
                    x,y->rcv(y,UCT!.exponent)));
orbits:=OrbitsDomain(aut,coho,gens,autmat);
return List(orbits,x->List(x,y->List(
                                 y,z->Int(z))));
end;

#& end #&


#& exponent #&

ExponentOfPFC:=function(G,R,K,UCT,hom)
local 	H,
	R2,
	K2,
	UCT2,
	m,
	phi,
	chainmap,
	res,
	list;

list:=[];
for H in Filtered(AllSubgroups(G),
		IsCyclic and IsNonTrivial) do
R2:=ResolutionFiniteGroup(H,4);
K2:=TensorWithIntegers(R2);
UCT2:=UniversalCoefficientsTheorem(K2,3);
m:=UCT2!.exponent;
phi:=GroupHomomorphismByFunction(H,G,h->h);
chainmap:=TensorWithIntegers(
	EquivariantChainMap(R2,R,phi));
res:=GroupCohomologyFunctor(
      K2,UCT!.complex,chainmap,3)!.mapping(
		hom);
Add(list,[Size(H),CocycleOrder(res,m)]);
od;
return Lcm(List(list,x->x[1]*x[2]));
end;

#& end #&




#& pi #&

ValuesOfPiSymbols:=function(G,R,K,UCT,x,f)
local n,e,o,list1,list2,m,i,j;
n:=UCT!.exponent;
e:=Exponent(G);
o:=Order(x);
list1:=[1];
for m in [1..o] do
Add(list2,
 factors[Size(list2)]*E(n)^f(x,x^(m-1),x));
od;
list2:=[1];
for j in [1..e] do
Add(list2,(list1[o+1]^QuoInt(j,o))
          *(list1[RemInt(j,o)+1]));
od;
return list2;
end;


#& end #&

#& twist #&

Twist:=function (class, C, lift, chi)
	return lift( class ) ^ chi / 
		lift( class ^ Size( C ) ) ^ chi;
end;



#& end #&


alphag:=function(g)
local f,x,y;
return function(f)
	return function(x,y)
    		return f(x,y,g)-f(x,y*g*y^(-1),y)+f(x*y*g*y^(-1)*x^(-1),x,y);end;
	end;
end;

#& indicator #&

FSIforDoubles:=function(
	G, 
	f, 
	class, 
	proj, i, 
	n, 
	listG, 
	pivalues, 
	m
	)
local e,C,Gaa,lift,chi,sum,x;
e:=Exponent(G);
C:=Centralizer(G,class);
Gaa:=proj!.Gaa;
lift:=proj!.lift;
chi:=proj!.table[i];
sum:=0;
for x in G do    
	if (class*x)^m=x^m then 
sum:=sum + E(n)^alphag(x^m)(f)(class,x)*
((pivalues[Position(listG,class*x)][e+1]/
pivalues[Position(listG,x)][e+1])^QuoInt(m,e))*
(pivalues[Position(listG,class*x)][RemInt(m,e)+1]/
pivalues[Position(listG,x)][RemInt(m,e)+1])*
lift((x^m))^chi;
        fi;
od;
return sum/Size(C);
end;


#& end #&





GaugeInvariants:=function(G)
    local 	R,
		K,
		UCT,
		n,
		hom,
		lift,
		orbits,
		representatives,
		o,
		listG,
		list,
		j,
		h,
		omega,
		f,
		exp,
		conj,
		PC,
		pivalues,
		FSlist,
		l,
		indicators,
		twist,
		columns;
	R:=ResolutionFiniteGroup(G,4);
	K:=TensorWithIntegers(R);
    	UCT:=UniversalCoefficientsTheorem(K,3);
    	Print("Cohomology done\n\n");
    	n:=UCT!.exponent;
	hombasis:=UCT!.hombasis;
    	lift:=UCT!.lift;
    	orbits:=AutomorphismOrbits(G,R,K,UCT);
    	representatives:=List([1..Size(orbits)],i->(lift*orbits[i][1])mod n);
    	Print("automorphism orbits done\n\n");
	listG:=Enumerator(G);
    	list:=[];
    	for j in [1..Size(representatives)] do
		h:=orbits[j][1];
        	omega:=representatives[j];
        	f:=StandardCocycle(R,omega,3,n);
        	exp:=ExponentOfPFC(G,R,K,UCT,h);
		conj:=List(ConjugacyClasses(G),x->Representative(x));
        	PC:=List(conj,g->ProjectiveCharacters(Centralizer(G,g),alphag(g)(f),n));
        	pivalues:=List(listG,x->ValuesOfPiSymbols(G,R,K,UCT,x,f));
		FSlist:=DivisorsInt(exp);l:=Size(FSlist);FSlist:=Concatenation([0],FSlist{[1..l-1]});
        	indicators:=Concatenation(List([1..Size(PC)],i->List([1..Size(PC[i]!.table)],j->List(FSlist,m->FSIforDoubles(G,f,conj[i],PC[i],j,n,listG,pivalues,m)))));
        	twist:=Concatenation(List([1..Size(PC)],i->List([1..Size(PC[i]!.table)],j->Twist(conj[i],Centralizer(G,conj[i]),PC[i]!.lift,PC[i]!.table[j]))));
        	SortParallel(indicators,twist);
        	columns:=Size(FSlist);
        	Add(list,[representatives[j],FSlist,Collected(List([1..Size(indicators)],i->Concatenation(List([1..columns],j->indicators[i][j]),[twist[i]])))]);
		Print(String(representatives[j]),"done\n\n");
    	od;
    	return rec( UCT:=UCT, orbits:=orbits, representatives:=representatives, Table:=list);
end;






#& ptindicators #&
		

FSIforPFC:=function(G,omega)
local R,K,UCT,n,f,listG,pivalues,list,result,nu,g,m;
R:=ResolutionFiniteGroup(G,4);
K:=TensorWithIntegers(R);
UCT:=UniversalCoefficientsTheorem(K,3);
n:=UCT!.exponent;
f:=StandardCocycle(R,omega,3,n);
listG:=Enumerator(G);
pivalues:=List(listG,x->
               ValuesOfPiSymbols(G,R,K,UCT,x,f));
nu:=function(G,f,x,m)
local o;
if x^m=One(G) then
  o:=Order(x);
  return (pivalues[Position(listG,x)][o+1]
           ^(QuoInt(-m,o)-1))
    *pivalues[Position(listG,x)][o-RemInt(-m,o)+1];
else return 0;
fi;
end;
result:=[];
for g in G do
  list:=[];
  for m in [0..ExponentOfPFC(G,R,K,UCT,
            (UCT!.cyclebasis*omega) mod n)-1] do
    Add(list,nu(G,f,g,m));
  od;
  Add(result,list);
od;
return result;
end;

#& end #&

#& adaptation #&

d:=function(eta)
return
  function(g,h,k) 
   return eta(h,k)-eta(g*h,k)+eta(g,h*k)-eta(g,h);
  end;
end;

moduloQ:=function(G,N,Q,g)
local p;
for p in Q do
 if p^(-1)*g in N then return [p,p^(-1)*g]; fi;
od;
end;

IsAdaptatedCocycle:=function(G,R,N,omega)
    local e,g,h,k,bool;
    e:=R!.exponent;
    bool:=true;
    for g in G do
        for h in G do
            for k in N do
                bool:=(omega(g,h,k)=0 mod e);
                if bool=false 
                then 
                    break;
                fi;
            od;
        od;
    od;
    return bool;
end;



FirstStep:=function(G,R,N,omega)
local Q,eta_1,omega_0;
Q:=List(RightCosets(G,N),x->Representative(x));
eta_1:=function(g,h)
 return omega(moduloQ(G,N,Q,g)[1],
         moduloQ(G,N,Q,g)[2],moduloQ(G,N,Q,h)[2]);
end;
omega_0:=function(g,h,k) 
 return omega(g,h,k)+d(eta_1)(g,h,k); 
end;
return omega_0;
end;

SecondStep:=function(G,R,N,omega_0)
local Q,eta_2,omega_1;
Q:=List(RightCosets(G,N),x->Representative(x));
eta_2:=function(g,h) 
 return 
  -(omega_0(g,moduloQ(G,N,Q,h)[1],
            moduloQ(G,N,Q,h)[2]))
  +omega_0(moduloQ(G,N,Q,g)[1],moduloQ(G,N,Q,g)[2],
   moduloQ(G,N,Q,h)[1]);
 end;
omega_1:=function(g,h,k) 
 return omega_0(g,h,k)+d(eta_2)(g,h,k);
end;
return omega_1;
end;

#& end #&




#& gtindicators #&

FSIforAGT:=function(G,N,omega)
local R,K,UCT,hom,FSexp,f,omega_0,omega_1,
      cosets,PC,list,result,i,j,nu,m;
R:=ResolutionFiniteGroup(G,4);
K:=TensorWithIntegers(R);
UCT:=UniversalCoefficientsTheorem(K,3);
hom:=(UCT!.cyclebasis*omega) mod UCT!.exponent;
FSexp:=ExponentOfPFC(G,R,K,UCT,hom);
f:=StandardCocycle(R,omega,3,UCT!.exponent);
omega_0:=FirstStep(G,R,N,f);
omega_1:=SecondStep(G,R,N,omega_0);
cosets:=DoubleCosets(G,N,N);
PC:=List(cosets,g->ProjectiveCharacters(
  Intersection(N,N^(Representative(g)^(-1))),
  function(x,y)
   return omega_1(x,y,Representative(g));
  end,
  UCT!.exponent));
nu:=function(G,UCT,N,f,class,list,i,m)
 local x,sum,S,Pf,lift,chi,listG,pivalues,r,o;
 S:=Intersection(N,
     N^(Representative(class)^(-1)));
 Gaa:=list!.Gaa;
 lift:=list!.lift;
 chi:=list!.table[i];
 listG:-Enumerator(G);
 pivalues:=List(listG,x->
            ValuesOfPiSymbols(G,R,K,UCT,x,f));
 sum:=0;
 for r in N do
  if (Representative(class)*r)^m in S 
  then 
   o:=Order(Representative(class)*r);
   sum:=sum 
    +(pivalues[Position(listG,
                Representative(class)*r)][o+1]
      ^(QuoInt(-m,o)-1))
    *pivalues[Position(listG,
                       Representative(class)*r)
      [o+RemInt(-m,o)+1]
    *lift(((Representative(class)*r)^(-m)))^chi; 
  fi;
 od;
return sum/Size(S);
end;
result:=[];
for i  in [1..Size(PC)]  do
 for j in [1..Size(PC[i]!.table)] do
  list:=[];
  for m in [0..FSexp-1] do
   Add(list,
       nu(G,UCT,N,omega_1,cosets[i],PC[i],j,m));
  od;
  Add(result,list);
 od;
od;
return result;
end;
#& end #&

\end{lstlisting}\label{lst:UCTcode}
\begin{remark}
There is already a command in HAP computing the torsion of the cohomology (or homology) of a finite group with values in $\mathbb{Z}$. Our command (whose primary result is of course a system of generators and relations) always gives the torsion in the form of the invariant factors of the cohomology group, due to the use of the Smith Normal Form; the HAP command sometimes gives the primary factors instead.
\end{remark}
The generators for $H^n(K,A)$ are represented as vectors; the order of the cohomology class corresponding to such a vector is given by the function:
\begin{lstlisting}[linerange=order-end]
LoadPackage("hapco");


InstallMethod( OneMutable,true,[IsCcElement],0,c -> One(InCcGroup(c)) );




ComplementIntMatWithTransforms:=function( full,sub )
local F, S, M, r, T, c, R;
  F := BaseIntMat( full );
  if IsEmpty( sub ) or IsZero( sub ) then
    return rec( complement := F, sub := [], moduli := [] );
  fi;
  S := BaseIntersectionIntMats( F, sub );
  if S <> BaseIntMat( sub ) then
    Error( "sub must be submodule of full" );
  fi;
  # find T such that BaseIntMat(T*F) = S
  M := Concatenation( F, S );
  T := NormalFormIntMat( M, 4 ).rowtrans^-1;
  T := T{[Length(F)+1..Length(T)]}{[1..Length(F)]};
  c := NormalFormIntMat(T*F,6)!.rowtrans;

  # r.rowtrans * T * F = r.normal * r.coltrans^-1 * F
  r := NormalFormIntMat( T, 13 );
  M := r.coltrans^-1 * F;
  R := rec( complement := BaseIntMat( M{[1+r.rank..Length(M)]} ),
            sub := r.rowtrans * T * F,
            moduli := List( [1..r.rank], i -> r.normal[i][i] ) ,
	    fulltrans := r.coltrans^(-1),
	    subtrans := r!.rowtrans*c^(-1));
  return R;
end;


#& cohomology #&

UniversalCoefficientsTheorem:=function(K,n)
local ZZ,     # a basis for Z_n(K)
      B,      # a basis for B_n(K)
      r,      # the ABT for Z_n(K) and B_n(K)
      C,      # the matrix C
      D,      # the list of the invariant factors
      m,      # the exponent of the H_n(K)
      homlist,# a "basis" for Hom(H_n(K),Z/mZ)
      r2, R1,C1,
              # the different variables implicated
              # in the ABT for K_n and Z_n(K)
		
      I,      # the matrix I
      cocyclelist;
              # a "basis" for H^n(K,C^*)


	
ZZ:=NullspaceIntMat(TransposedMat(
                BoundaryMatrix(K,n)));
if not Size(ZZ)>0 then return 
	rec(	complex:=K,
		cyclebasis:=ZZ,
		cycletrans:=[],
		hombasis := [],
		cocyclebasis := [],
		torsion := [1],
		exponent := 1,
		lift:=[]
	);
fi;
B:=BaseIntMat(TransposedMat(
                 BoundaryMatrix(K,n+1)));
r:=ComplementIntMatWithTransforms(ZZ,B);
C:=r!.fulltrans;
D:=r!.moduli;
m:=D[Size(D)];
homlist:=[];
List(Filtered(D,x->not x=1),x->
	ListWithIdenticalEntries(Size(ZZ),0));
for i in [1..Size(D)] do
	if D[i]=1 then continue;fi;
	hom:=ListWithIdenticalEntries(Size(ZZ),0);
	hom[i]:=m/D[i];
	Add(homlist,hom);
od;
homlist:=List(homlist,x->((C^(-1))*x) mod m);

r2:=ComplementIntMatWithTransforms(
              IdentityMat(K!.dimension(3)),ZZ);
C1:=r2!.subtrans;
R1:=r2!.fulltrans;

I:=Concatenation(IdentityMat(Size(ZZ)),
     NullMat(K!.dimension(n)-Size(ZZ),Size(ZZ)));
cocyclelist:=List(homlist,x->(R1^(-1)*(I*(C1*x))) 
                     mod m);
	
return( rec(
	complex:=K,
	cyclebasis:=ZZ,
	cycletrans:=C,
	hombasis := homlist,
	cocyclebasis := cocyclelist,
	torsion := Filtered(D,x->not x=1),
	exponent := m,
	lift:=R1^(-1)*I*C1
		)
	);
end;


#& end #&


#& order #&

CocycleOrder:=function(f,m)
	local list;
	list:=List(Filtered(f,x->not x=0),
                   x->m/Gcd(x,m));
	if list=[] 
		then return 1;
		else return Maximum(list);
	fi;
end;    

#& end #&   

#& functoriality #&

GroupCohomologyFunctor:=function(K1,K2,phi,n)
local 	UCT1, 	# UCT for K1
	m1, 	# the exponent of H_n(K1)
	Z1,	# the adapted basis of Z_n(K1)
	UCT2, 	# UCT for K2
	m2, 	# the exponent of H_n(K2)
	Z2,	# the adapted basis of Z_n(K2)
	p, 	# p=lcm(m1,m2)
	F, 	# the matrix F
	hphi; 	# the morphism H_n(phi)
		

UCT1:=UniversalCoefficientsTheorem(K1,n);;
m1:=UCT1!.exponent;;
Z1:=UCT1!.cyclebasis;
s:=NormalFormIntMat(Z1,13);
UCT2:=UniversalCoefficientsTheorem(K2,n);;
m2:=UCT2!.exponent;;
Z2:=UCT2!.cyclebasis;
F:=List(Z1,
	x->SolutionIntMat(Z2,phi!.mapping(x,n)));;
p:=Lcm(m1,m2);;

hphi:=function(f)
	return ((F*(m1*f))/m2) mod m1;
end;
return rec(
	matrix:=F,
	mapping:=hphi,
	p:=p,
	m1:=m1,
	m2:=m2);
end;

#& end #&


#& proj #&

ProjectiveCharacters:=function(G,alpha,o)
local A,r,Ga,Gaa,ire,t,lift,testpro;

# we represent the cyclic group as a G-Outer group

A:=GOuterGroup();
SetActingGroup(A,G);
if o=1 then
    SetActedGroup(A,GroupByGenerators(
	[One(CyclicGroup(o))]));
    else
        SetActedGroup(A,CyclicGroup(o));
fi;
SetOuterAction(A,function(g,x) return x; end);

# we represent the cocycle f as a 
# standard 2-cocycle in HAP 

r:=Standard2Cocycle();
SetActingGroup(r,A!.ActingGroup);
SetMapping(r,function(x,y) return 
	(A!.ActedGroup.1)^alpha(x,y);
	end);
SetModule(r,A);


# now, we construct the extension of G by A along r
# we also give its character table

Ga:=CcGroup(A,r);
Gaa:=UnderlyingGroup(Ga);
ire:=Irr(Gaa);

# we give the element in Gaa that corresponds to 
# the generators of the cyclic group, and a 
# section of the epi of Gaa on G


t:=CcElement(FamilyObj(One(Ga)),GeneratorsOfGroup(
   Ga!.Fibre)[1],One(Ga!.Base),InCcGroup(One(Ga)));
lift:=function(g)
      return CcElement(FamilyObj(One(Ga)),
	One(Ga!.Fibre),g,InCcGroup(One(Ga)));
end;

# finally, we give the function that tests for a 
# charater of Gaa if it comes from a projective 
# character of G

testpro:=function(l)
	local list,chi;
	list:=[];
	for chi in l do
      		if t^chi/E(o)=(t^Size(Ga))^chi 
			then Add(list,chi); 
		fi;
	od;
	return list;
end;

return rec(
	Gaa:=Gaa,
	lift:=lift,
	table:=testpro(ire)
	);
end;

#& end #&


#& automorphisms #&

AutomorphismOrbits:=function(G,R,K,UCT)
local 	rcv,
	coho,
	aut,
	gens,
	autmat,
	orbits;

rcv:=function(x,n) 
     return List(x,y->ZmodnZObj(y,n));
     end;
coho:=NearAdditiveGroup(
      List(UCT!.hombasis,x->rcv(x,UCT!.exponent)));
aut:=AutomorphismGroup(G);
gens:=GeneratorsOfGroup(aut);
autmat:=List(
         gens,phi->TransposedMat(
          GroupCohomologyFunctor(
           K,K,TensorWithIntegers(
	    EquivariantChainMap(R,R,phi)),
           3)!.matrix));
autmat:=List(
         autmat,x->List(
                    x,y->rcv(y,UCT!.exponent)));
orbits:=OrbitsDomain(aut,coho,gens,autmat);
return List(orbits,x->List(x,y->List(
                                 y,z->Int(z))));
end;

#& end #&


#& exponent #&

ExponentOfPFC:=function(G,R,K,UCT,hom)
local 	H,
	R2,
	K2,
	UCT2,
	m,
	phi,
	chainmap,
	res,
	list;

list:=[];
for H in Filtered(AllSubgroups(G),
		IsCyclic and IsNonTrivial) do
R2:=ResolutionFiniteGroup(H,4);
K2:=TensorWithIntegers(R2);
UCT2:=UniversalCoefficientsTheorem(K2,3);
m:=UCT2!.exponent;
phi:=GroupHomomorphismByFunction(H,G,h->h);
chainmap:=TensorWithIntegers(
	EquivariantChainMap(R2,R,phi));
res:=GroupCohomologyFunctor(
      K2,UCT!.complex,chainmap,3)!.mapping(
		hom);
Add(list,[Size(H),CocycleOrder(res,m)]);
od;
return Lcm(List(list,x->x[1]*x[2]));
end;

#& end #&




#& pi #&

ValuesOfPiSymbols:=function(G,R,K,UCT,x,f)
local n,e,o,list1,list2,m,i,j;
n:=UCT!.exponent;
e:=Exponent(G);
o:=Order(x);
list1:=[1];
for m in [1..o] do
Add(list2,
 factors[Size(list2)]*E(n)^f(x,x^(m-1),x));
od;
list2:=[1];
for j in [1..e] do
Add(list2,(list1[o+1]^QuoInt(j,o))
          *(list1[RemInt(j,o)+1]));
od;
return list2;
end;


#& end #&

#& twist #&

Twist:=function (class, C, lift, chi)
	return lift( class ) ^ chi / 
		lift( class ^ Size( C ) ) ^ chi;
end;



#& end #&


alphag:=function(g)
local f,x,y;
return function(f)
	return function(x,y)
    		return f(x,y,g)-f(x,y*g*y^(-1),y)+f(x*y*g*y^(-1)*x^(-1),x,y);end;
	end;
end;

#& indicator #&

FSIforDoubles:=function(
	G, 
	f, 
	class, 
	proj, i, 
	n, 
	listG, 
	pivalues, 
	m
	)
local e,C,Gaa,lift,chi,sum,x;
e:=Exponent(G);
C:=Centralizer(G,class);
Gaa:=proj!.Gaa;
lift:=proj!.lift;
chi:=proj!.table[i];
sum:=0;
for x in G do    
	if (class*x)^m=x^m then 
sum:=sum + E(n)^alphag(x^m)(f)(class,x)*
((pivalues[Position(listG,class*x)][e+1]/
pivalues[Position(listG,x)][e+1])^QuoInt(m,e))*
(pivalues[Position(listG,class*x)][RemInt(m,e)+1]/
pivalues[Position(listG,x)][RemInt(m,e)+1])*
lift((x^m))^chi;
        fi;
od;
return sum/Size(C);
end;


#& end #&





GaugeInvariants:=function(G)
    local 	R,
		K,
		UCT,
		n,
		hom,
		lift,
		orbits,
		representatives,
		o,
		listG,
		list,
		j,
		h,
		omega,
		f,
		exp,
		conj,
		PC,
		pivalues,
		FSlist,
		l,
		indicators,
		twist,
		columns;
	R:=ResolutionFiniteGroup(G,4);
	K:=TensorWithIntegers(R);
    	UCT:=UniversalCoefficientsTheorem(K,3);
    	Print("Cohomology done\n\n");
    	n:=UCT!.exponent;
	hombasis:=UCT!.hombasis;
    	lift:=UCT!.lift;
    	orbits:=AutomorphismOrbits(G,R,K,UCT);
    	representatives:=List([1..Size(orbits)],i->(lift*orbits[i][1])mod n);
    	Print("automorphism orbits done\n\n");
	listG:=Enumerator(G);
    	list:=[];
    	for j in [1..Size(representatives)] do
		h:=orbits[j][1];
        	omega:=representatives[j];
        	f:=StandardCocycle(R,omega,3,n);
        	exp:=ExponentOfPFC(G,R,K,UCT,h);
		conj:=List(ConjugacyClasses(G),x->Representative(x));
        	PC:=List(conj,g->ProjectiveCharacters(Centralizer(G,g),alphag(g)(f),n));
        	pivalues:=List(listG,x->ValuesOfPiSymbols(G,R,K,UCT,x,f));
		FSlist:=DivisorsInt(exp);l:=Size(FSlist);FSlist:=Concatenation([0],FSlist{[1..l-1]});
        	indicators:=Concatenation(List([1..Size(PC)],i->List([1..Size(PC[i]!.table)],j->List(FSlist,m->FSIforDoubles(G,f,conj[i],PC[i],j,n,listG,pivalues,m)))));
        	twist:=Concatenation(List([1..Size(PC)],i->List([1..Size(PC[i]!.table)],j->Twist(conj[i],Centralizer(G,conj[i]),PC[i]!.lift,PC[i]!.table[j]))));
        	SortParallel(indicators,twist);
        	columns:=Size(FSlist);
        	Add(list,[representatives[j],FSlist,Collected(List([1..Size(indicators)],i->Concatenation(List([1..columns],j->indicators[i][j]),[twist[i]])))]);
		Print(String(representatives[j]),"done\n\n");
    	od;
    	return rec( UCT:=UCT, orbits:=orbits, representatives:=representatives, Table:=list);
end;






#& ptindicators #&
		

FSIforPFC:=function(G,omega)
local R,K,UCT,n,f,listG,pivalues,list,result,nu,g,m;
R:=ResolutionFiniteGroup(G,4);
K:=TensorWithIntegers(R);
UCT:=UniversalCoefficientsTheorem(K,3);
n:=UCT!.exponent;
f:=StandardCocycle(R,omega,3,n);
listG:=Enumerator(G);
pivalues:=List(listG,x->
               ValuesOfPiSymbols(G,R,K,UCT,x,f));
nu:=function(G,f,x,m)
local o;
if x^m=One(G) then
  o:=Order(x);
  return (pivalues[Position(listG,x)][o+1]
           ^(QuoInt(-m,o)-1))
    *pivalues[Position(listG,x)][o-RemInt(-m,o)+1];
else return 0;
fi;
end;
result:=[];
for g in G do
  list:=[];
  for m in [0..ExponentOfPFC(G,R,K,UCT,
            (UCT!.cyclebasis*omega) mod n)-1] do
    Add(list,nu(G,f,g,m));
  od;
  Add(result,list);
od;
return result;
end;

#& end #&

#& adaptation #&

d:=function(eta)
return
  function(g,h,k) 
   return eta(h,k)-eta(g*h,k)+eta(g,h*k)-eta(g,h);
  end;
end;

moduloQ:=function(G,N,Q,g)
local p;
for p in Q do
 if p^(-1)*g in N then return [p,p^(-1)*g]; fi;
od;
end;

IsAdaptatedCocycle:=function(G,R,N,omega)
    local e,g,h,k,bool;
    e:=R!.exponent;
    bool:=true;
    for g in G do
        for h in G do
            for k in N do
                bool:=(omega(g,h,k)=0 mod e);
                if bool=false 
                then 
                    break;
                fi;
            od;
        od;
    od;
    return bool;
end;



FirstStep:=function(G,R,N,omega)
local Q,eta_1,omega_0;
Q:=List(RightCosets(G,N),x->Representative(x));
eta_1:=function(g,h)
 return omega(moduloQ(G,N,Q,g)[1],
         moduloQ(G,N,Q,g)[2],moduloQ(G,N,Q,h)[2]);
end;
omega_0:=function(g,h,k) 
 return omega(g,h,k)+d(eta_1)(g,h,k); 
end;
return omega_0;
end;

SecondStep:=function(G,R,N,omega_0)
local Q,eta_2,omega_1;
Q:=List(RightCosets(G,N),x->Representative(x));
eta_2:=function(g,h) 
 return 
  -(omega_0(g,moduloQ(G,N,Q,h)[1],
            moduloQ(G,N,Q,h)[2]))
  +omega_0(moduloQ(G,N,Q,g)[1],moduloQ(G,N,Q,g)[2],
   moduloQ(G,N,Q,h)[1]);
 end;
omega_1:=function(g,h,k) 
 return omega_0(g,h,k)+d(eta_2)(g,h,k);
end;
return omega_1;
end;

#& end #&




#& gtindicators #&

FSIforAGT:=function(G,N,omega)
local R,K,UCT,hom,FSexp,f,omega_0,omega_1,
      cosets,PC,list,result,i,j,nu,m;
R:=ResolutionFiniteGroup(G,4);
K:=TensorWithIntegers(R);
UCT:=UniversalCoefficientsTheorem(K,3);
hom:=(UCT!.cyclebasis*omega) mod UCT!.exponent;
FSexp:=ExponentOfPFC(G,R,K,UCT,hom);
f:=StandardCocycle(R,omega,3,UCT!.exponent);
omega_0:=FirstStep(G,R,N,f);
omega_1:=SecondStep(G,R,N,omega_0);
cosets:=DoubleCosets(G,N,N);
PC:=List(cosets,g->ProjectiveCharacters(
  Intersection(N,N^(Representative(g)^(-1))),
  function(x,y)
   return omega_1(x,y,Representative(g));
  end,
  UCT!.exponent));
nu:=function(G,UCT,N,f,class,list,i,m)
 local x,sum,S,Pf,lift,chi,listG,pivalues,r,o;
 S:=Intersection(N,
     N^(Representative(class)^(-1)));
 Gaa:=list!.Gaa;
 lift:=list!.lift;
 chi:=list!.table[i];
 listG:-Enumerator(G);
 pivalues:=List(listG,x->
            ValuesOfPiSymbols(G,R,K,UCT,x,f));
 sum:=0;
 for r in N do
  if (Representative(class)*r)^m in S 
  then 
   o:=Order(Representative(class)*r);
   sum:=sum 
    +(pivalues[Position(listG,
                Representative(class)*r)][o+1]
      ^(QuoInt(-m,o)-1))
    *pivalues[Position(listG,
                       Representative(class)*r)
      [o+RemInt(-m,o)+1]
    *lift(((Representative(class)*r)^(-m)))^chi; 
  fi;
 od;
return sum/Size(S);
end;
result:=[];
for i  in [1..Size(PC)]  do
 for j in [1..Size(PC[i]!.table)] do
  list:=[];
  for m in [0..FSexp-1] do
   Add(list,
       nu(G,UCT,N,omega_1,cosets[i],PC[i],j,m));
  od;
  Add(result,list);
 od;
od;
return result;
end;
#& end #&

\end{lstlisting}\label{lst:cocycleordercode}

Now, for a chain map $f$ from a complex $K_1$ to a complex $K_2$, and an integer $n$, we give the code that provides the matrix $F$ representing the homomorphism $H_n(f)$.
\begin{lstlisting}[linerange=functoriality-end]
LoadPackage("hapco");


InstallMethod( OneMutable,true,[IsCcElement],0,c -> One(InCcGroup(c)) );




ComplementIntMatWithTransforms:=function( full,sub )
local F, S, M, r, T, c, R;
  F := BaseIntMat( full );
  if IsEmpty( sub ) or IsZero( sub ) then
    return rec( complement := F, sub := [], moduli := [] );
  fi;
  S := BaseIntersectionIntMats( F, sub );
  if S <> BaseIntMat( sub ) then
    Error( "sub must be submodule of full" );
  fi;
  # find T such that BaseIntMat(T*F) = S
  M := Concatenation( F, S );
  T := NormalFormIntMat( M, 4 ).rowtrans^-1;
  T := T{[Length(F)+1..Length(T)]}{[1..Length(F)]};
  c := NormalFormIntMat(T*F,6)!.rowtrans;

  # r.rowtrans * T * F = r.normal * r.coltrans^-1 * F
  r := NormalFormIntMat( T, 13 );
  M := r.coltrans^-1 * F;
  R := rec( complement := BaseIntMat( M{[1+r.rank..Length(M)]} ),
            sub := r.rowtrans * T * F,
            moduli := List( [1..r.rank], i -> r.normal[i][i] ) ,
	    fulltrans := r.coltrans^(-1),
	    subtrans := r!.rowtrans*c^(-1));
  return R;
end;


#& cohomology #&

UniversalCoefficientsTheorem:=function(K,n)
local ZZ,     # a basis for Z_n(K)
      B,      # a basis for B_n(K)
      r,      # the ABT for Z_n(K) and B_n(K)
      C,      # the matrix C
      D,      # the list of the invariant factors
      m,      # the exponent of the H_n(K)
      homlist,# a "basis" for Hom(H_n(K),Z/mZ)
      r2, R1,C1,
              # the different variables implicated
              # in the ABT for K_n and Z_n(K)
		
      I,      # the matrix I
      cocyclelist;
              # a "basis" for H^n(K,C^*)


	
ZZ:=NullspaceIntMat(TransposedMat(
                BoundaryMatrix(K,n)));
if not Size(ZZ)>0 then return 
	rec(	complex:=K,
		cyclebasis:=ZZ,
		cycletrans:=[],
		hombasis := [],
		cocyclebasis := [],
		torsion := [1],
		exponent := 1,
		lift:=[]
	);
fi;
B:=BaseIntMat(TransposedMat(
                 BoundaryMatrix(K,n+1)));
r:=ComplementIntMatWithTransforms(ZZ,B);
C:=r!.fulltrans;
D:=r!.moduli;
m:=D[Size(D)];
homlist:=[];
List(Filtered(D,x->not x=1),x->
	ListWithIdenticalEntries(Size(ZZ),0));
for i in [1..Size(D)] do
	if D[i]=1 then continue;fi;
	hom:=ListWithIdenticalEntries(Size(ZZ),0);
	hom[i]:=m/D[i];
	Add(homlist,hom);
od;
homlist:=List(homlist,x->((C^(-1))*x) mod m);

r2:=ComplementIntMatWithTransforms(
              IdentityMat(K!.dimension(3)),ZZ);
C1:=r2!.subtrans;
R1:=r2!.fulltrans;

I:=Concatenation(IdentityMat(Size(ZZ)),
     NullMat(K!.dimension(n)-Size(ZZ),Size(ZZ)));
cocyclelist:=List(homlist,x->(R1^(-1)*(I*(C1*x))) 
                     mod m);
	
return( rec(
	complex:=K,
	cyclebasis:=ZZ,
	cycletrans:=C,
	hombasis := homlist,
	cocyclebasis := cocyclelist,
	torsion := Filtered(D,x->not x=1),
	exponent := m,
	lift:=R1^(-1)*I*C1
		)
	);
end;


#& end #&


#& order #&

CocycleOrder:=function(f,m)
	local list;
	list:=List(Filtered(f,x->not x=0),
                   x->m/Gcd(x,m));
	if list=[] 
		then return 1;
		else return Maximum(list);
	fi;
end;    

#& end #&   

#& functoriality #&

GroupCohomologyFunctor:=function(K1,K2,phi,n)
local 	UCT1, 	# UCT for K1
	m1, 	# the exponent of H_n(K1)
	Z1,	# the adapted basis of Z_n(K1)
	UCT2, 	# UCT for K2
	m2, 	# the exponent of H_n(K2)
	Z2,	# the adapted basis of Z_n(K2)
	p, 	# p=lcm(m1,m2)
	F, 	# the matrix F
	hphi; 	# the morphism H_n(phi)
		

UCT1:=UniversalCoefficientsTheorem(K1,n);;
m1:=UCT1!.exponent;;
Z1:=UCT1!.cyclebasis;
s:=NormalFormIntMat(Z1,13);
UCT2:=UniversalCoefficientsTheorem(K2,n);;
m2:=UCT2!.exponent;;
Z2:=UCT2!.cyclebasis;
F:=List(Z1,
	x->SolutionIntMat(Z2,phi!.mapping(x,n)));;
p:=Lcm(m1,m2);;

hphi:=function(f)
	return ((F*(m1*f))/m2) mod m1;
end;
return rec(
	matrix:=F,
	mapping:=hphi,
	p:=p,
	m1:=m1,
	m2:=m2);
end;

#& end #&


#& proj #&

ProjectiveCharacters:=function(G,alpha,o)
local A,r,Ga,Gaa,ire,t,lift,testpro;

# we represent the cyclic group as a G-Outer group

A:=GOuterGroup();
SetActingGroup(A,G);
if o=1 then
    SetActedGroup(A,GroupByGenerators(
	[One(CyclicGroup(o))]));
    else
        SetActedGroup(A,CyclicGroup(o));
fi;
SetOuterAction(A,function(g,x) return x; end);

# we represent the cocycle f as a 
# standard 2-cocycle in HAP 

r:=Standard2Cocycle();
SetActingGroup(r,A!.ActingGroup);
SetMapping(r,function(x,y) return 
	(A!.ActedGroup.1)^alpha(x,y);
	end);
SetModule(r,A);


# now, we construct the extension of G by A along r
# we also give its character table

Ga:=CcGroup(A,r);
Gaa:=UnderlyingGroup(Ga);
ire:=Irr(Gaa);

# we give the element in Gaa that corresponds to 
# the generators of the cyclic group, and a 
# section of the epi of Gaa on G


t:=CcElement(FamilyObj(One(Ga)),GeneratorsOfGroup(
   Ga!.Fibre)[1],One(Ga!.Base),InCcGroup(One(Ga)));
lift:=function(g)
      return CcElement(FamilyObj(One(Ga)),
	One(Ga!.Fibre),g,InCcGroup(One(Ga)));
end;

# finally, we give the function that tests for a 
# charater of Gaa if it comes from a projective 
# character of G

testpro:=function(l)
	local list,chi;
	list:=[];
	for chi in l do
      		if t^chi/E(o)=(t^Size(Ga))^chi 
			then Add(list,chi); 
		fi;
	od;
	return list;
end;

return rec(
	Gaa:=Gaa,
	lift:=lift,
	table:=testpro(ire)
	);
end;

#& end #&


#& automorphisms #&

AutomorphismOrbits:=function(G,R,K,UCT)
local 	rcv,
	coho,
	aut,
	gens,
	autmat,
	orbits;

rcv:=function(x,n) 
     return List(x,y->ZmodnZObj(y,n));
     end;
coho:=NearAdditiveGroup(
      List(UCT!.hombasis,x->rcv(x,UCT!.exponent)));
aut:=AutomorphismGroup(G);
gens:=GeneratorsOfGroup(aut);
autmat:=List(
         gens,phi->TransposedMat(
          GroupCohomologyFunctor(
           K,K,TensorWithIntegers(
	    EquivariantChainMap(R,R,phi)),
           3)!.matrix));
autmat:=List(
         autmat,x->List(
                    x,y->rcv(y,UCT!.exponent)));
orbits:=OrbitsDomain(aut,coho,gens,autmat);
return List(orbits,x->List(x,y->List(
                                 y,z->Int(z))));
end;

#& end #&


#& exponent #&

ExponentOfPFC:=function(G,R,K,UCT,hom)
local 	H,
	R2,
	K2,
	UCT2,
	m,
	phi,
	chainmap,
	res,
	list;

list:=[];
for H in Filtered(AllSubgroups(G),
		IsCyclic and IsNonTrivial) do
R2:=ResolutionFiniteGroup(H,4);
K2:=TensorWithIntegers(R2);
UCT2:=UniversalCoefficientsTheorem(K2,3);
m:=UCT2!.exponent;
phi:=GroupHomomorphismByFunction(H,G,h->h);
chainmap:=TensorWithIntegers(
	EquivariantChainMap(R2,R,phi));
res:=GroupCohomologyFunctor(
      K2,UCT!.complex,chainmap,3)!.mapping(
		hom);
Add(list,[Size(H),CocycleOrder(res,m)]);
od;
return Lcm(List(list,x->x[1]*x[2]));
end;

#& end #&




#& pi #&

ValuesOfPiSymbols:=function(G,R,K,UCT,x,f)
local n,e,o,list1,list2,m,i,j;
n:=UCT!.exponent;
e:=Exponent(G);
o:=Order(x);
list1:=[1];
for m in [1..o] do
Add(list2,
 factors[Size(list2)]*E(n)^f(x,x^(m-1),x));
od;
list2:=[1];
for j in [1..e] do
Add(list2,(list1[o+1]^QuoInt(j,o))
          *(list1[RemInt(j,o)+1]));
od;
return list2;
end;


#& end #&

#& twist #&

Twist:=function (class, C, lift, chi)
	return lift( class ) ^ chi / 
		lift( class ^ Size( C ) ) ^ chi;
end;



#& end #&


alphag:=function(g)
local f,x,y;
return function(f)
	return function(x,y)
    		return f(x,y,g)-f(x,y*g*y^(-1),y)+f(x*y*g*y^(-1)*x^(-1),x,y);end;
	end;
end;

#& indicator #&

FSIforDoubles:=function(
	G, 
	f, 
	class, 
	proj, i, 
	n, 
	listG, 
	pivalues, 
	m
	)
local e,C,Gaa,lift,chi,sum,x;
e:=Exponent(G);
C:=Centralizer(G,class);
Gaa:=proj!.Gaa;
lift:=proj!.lift;
chi:=proj!.table[i];
sum:=0;
for x in G do    
	if (class*x)^m=x^m then 
sum:=sum + E(n)^alphag(x^m)(f)(class,x)*
((pivalues[Position(listG,class*x)][e+1]/
pivalues[Position(listG,x)][e+1])^QuoInt(m,e))*
(pivalues[Position(listG,class*x)][RemInt(m,e)+1]/
pivalues[Position(listG,x)][RemInt(m,e)+1])*
lift((x^m))^chi;
        fi;
od;
return sum/Size(C);
end;


#& end #&





GaugeInvariants:=function(G)
    local 	R,
		K,
		UCT,
		n,
		hom,
		lift,
		orbits,
		representatives,
		o,
		listG,
		list,
		j,
		h,
		omega,
		f,
		exp,
		conj,
		PC,
		pivalues,
		FSlist,
		l,
		indicators,
		twist,
		columns;
	R:=ResolutionFiniteGroup(G,4);
	K:=TensorWithIntegers(R);
    	UCT:=UniversalCoefficientsTheorem(K,3);
    	Print("Cohomology done\n\n");
    	n:=UCT!.exponent;
	hombasis:=UCT!.hombasis;
    	lift:=UCT!.lift;
    	orbits:=AutomorphismOrbits(G,R,K,UCT);
    	representatives:=List([1..Size(orbits)],i->(lift*orbits[i][1])mod n);
    	Print("automorphism orbits done\n\n");
	listG:=Enumerator(G);
    	list:=[];
    	for j in [1..Size(representatives)] do
		h:=orbits[j][1];
        	omega:=representatives[j];
        	f:=StandardCocycle(R,omega,3,n);
        	exp:=ExponentOfPFC(G,R,K,UCT,h);
		conj:=List(ConjugacyClasses(G),x->Representative(x));
        	PC:=List(conj,g->ProjectiveCharacters(Centralizer(G,g),alphag(g)(f),n));
        	pivalues:=List(listG,x->ValuesOfPiSymbols(G,R,K,UCT,x,f));
		FSlist:=DivisorsInt(exp);l:=Size(FSlist);FSlist:=Concatenation([0],FSlist{[1..l-1]});
        	indicators:=Concatenation(List([1..Size(PC)],i->List([1..Size(PC[i]!.table)],j->List(FSlist,m->FSIforDoubles(G,f,conj[i],PC[i],j,n,listG,pivalues,m)))));
        	twist:=Concatenation(List([1..Size(PC)],i->List([1..Size(PC[i]!.table)],j->Twist(conj[i],Centralizer(G,conj[i]),PC[i]!.lift,PC[i]!.table[j]))));
        	SortParallel(indicators,twist);
        	columns:=Size(FSlist);
        	Add(list,[representatives[j],FSlist,Collected(List([1..Size(indicators)],i->Concatenation(List([1..columns],j->indicators[i][j]),[twist[i]])))]);
		Print(String(representatives[j]),"done\n\n");
    	od;
    	return rec( UCT:=UCT, orbits:=orbits, representatives:=representatives, Table:=list);
end;






#& ptindicators #&
		

FSIforPFC:=function(G,omega)
local R,K,UCT,n,f,listG,pivalues,list,result,nu,g,m;
R:=ResolutionFiniteGroup(G,4);
K:=TensorWithIntegers(R);
UCT:=UniversalCoefficientsTheorem(K,3);
n:=UCT!.exponent;
f:=StandardCocycle(R,omega,3,n);
listG:=Enumerator(G);
pivalues:=List(listG,x->
               ValuesOfPiSymbols(G,R,K,UCT,x,f));
nu:=function(G,f,x,m)
local o;
if x^m=One(G) then
  o:=Order(x);
  return (pivalues[Position(listG,x)][o+1]
           ^(QuoInt(-m,o)-1))
    *pivalues[Position(listG,x)][o-RemInt(-m,o)+1];
else return 0;
fi;
end;
result:=[];
for g in G do
  list:=[];
  for m in [0..ExponentOfPFC(G,R,K,UCT,
            (UCT!.cyclebasis*omega) mod n)-1] do
    Add(list,nu(G,f,g,m));
  od;
  Add(result,list);
od;
return result;
end;

#& end #&

#& adaptation #&

d:=function(eta)
return
  function(g,h,k) 
   return eta(h,k)-eta(g*h,k)+eta(g,h*k)-eta(g,h);
  end;
end;

moduloQ:=function(G,N,Q,g)
local p;
for p in Q do
 if p^(-1)*g in N then return [p,p^(-1)*g]; fi;
od;
end;

IsAdaptatedCocycle:=function(G,R,N,omega)
    local e,g,h,k,bool;
    e:=R!.exponent;
    bool:=true;
    for g in G do
        for h in G do
            for k in N do
                bool:=(omega(g,h,k)=0 mod e);
                if bool=false 
                then 
                    break;
                fi;
            od;
        od;
    od;
    return bool;
end;



FirstStep:=function(G,R,N,omega)
local Q,eta_1,omega_0;
Q:=List(RightCosets(G,N),x->Representative(x));
eta_1:=function(g,h)
 return omega(moduloQ(G,N,Q,g)[1],
         moduloQ(G,N,Q,g)[2],moduloQ(G,N,Q,h)[2]);
end;
omega_0:=function(g,h,k) 
 return omega(g,h,k)+d(eta_1)(g,h,k); 
end;
return omega_0;
end;

SecondStep:=function(G,R,N,omega_0)
local Q,eta_2,omega_1;
Q:=List(RightCosets(G,N),x->Representative(x));
eta_2:=function(g,h) 
 return 
  -(omega_0(g,moduloQ(G,N,Q,h)[1],
            moduloQ(G,N,Q,h)[2]))
  +omega_0(moduloQ(G,N,Q,g)[1],moduloQ(G,N,Q,g)[2],
   moduloQ(G,N,Q,h)[1]);
 end;
omega_1:=function(g,h,k) 
 return omega_0(g,h,k)+d(eta_2)(g,h,k);
end;
return omega_1;
end;

#& end #&




#& gtindicators #&

FSIforAGT:=function(G,N,omega)
local R,K,UCT,hom,FSexp,f,omega_0,omega_1,
      cosets,PC,list,result,i,j,nu,m;
R:=ResolutionFiniteGroup(G,4);
K:=TensorWithIntegers(R);
UCT:=UniversalCoefficientsTheorem(K,3);
hom:=(UCT!.cyclebasis*omega) mod UCT!.exponent;
FSexp:=ExponentOfPFC(G,R,K,UCT,hom);
f:=StandardCocycle(R,omega,3,UCT!.exponent);
omega_0:=FirstStep(G,R,N,f);
omega_1:=SecondStep(G,R,N,omega_0);
cosets:=DoubleCosets(G,N,N);
PC:=List(cosets,g->ProjectiveCharacters(
  Intersection(N,N^(Representative(g)^(-1))),
  function(x,y)
   return omega_1(x,y,Representative(g));
  end,
  UCT!.exponent));
nu:=function(G,UCT,N,f,class,list,i,m)
 local x,sum,S,Pf,lift,chi,listG,pivalues,r,o;
 S:=Intersection(N,
     N^(Representative(class)^(-1)));
 Gaa:=list!.Gaa;
 lift:=list!.lift;
 chi:=list!.table[i];
 listG:-Enumerator(G);
 pivalues:=List(listG,x->
            ValuesOfPiSymbols(G,R,K,UCT,x,f));
 sum:=0;
 for r in N do
  if (Representative(class)*r)^m in S 
  then 
   o:=Order(Representative(class)*r);
   sum:=sum 
    +(pivalues[Position(listG,
                Representative(class)*r)][o+1]
      ^(QuoInt(-m,o)-1))
    *pivalues[Position(listG,
                       Representative(class)*r)
      [o+RemInt(-m,o)+1]
    *lift(((Representative(class)*r)^(-m)))^chi; 
  fi;
 od;
return sum/Size(S);
end;
result:=[];
for i  in [1..Size(PC)]  do
 for j in [1..Size(PC[i]!.table)] do
  list:=[];
  for m in [0..FSexp-1] do
   Add(list,
       nu(G,UCT,N,omega_1,cosets[i],PC[i],j,m));
  od;
  Add(result,list);
 od;
od;
return result;
end;
#& end #&

\end{lstlisting}\label{lst:functorcode}
\section{Frobenius-Schur indicators of group-theoretical categories}
\label{sec:context}

\subsection{Preliminaries}
\label{sec:gtc}

A \emph{fusion category} (over $\mathbb{C}$) is a rigid semi-simple, $\mathbb{C}$-linear and (non-strict) monoidal category with finitely many classes of isomorphism of simple objects. We also require that $ \Hom(s,t)=k$ if $s$ and $t$ are isomorphic simple objects, and $0$ otherwise. A \emph{pivotal structure} in a fusion category is a monoidal natural isomorphism $j: Id \rightarrow (\ )^{**}$.

Recall also that a \emph{braiding} in a monoidal category is a family of natural isomorphisms $c_{X,Y}:X\otimes Y \rightarrow Y\otimes X$ satisfying an hexagon axiom. A monoidal category equipped by a braiding is called \emph{braided}. A \emph{twist} in a braided category is a family of natural isomorphisms $\theta_V: V \rightarrow V^*$ such that $\theta_{X \otimes Y}=(\theta_X \otimes \theta_Y)c_{X,Y}c_{Y,X}$ and $\theta_{X^*}=(\theta_X)^*$. A braided category with twist is called a \emph{ribbon category}. For more background on braided and ribbon categories, we refer the reader to  \cite{Kassel}. 

The Drinfeld center $\mathcal{Z}(\C)$ of a monoidal category $\C$ is the category whose objects are couple $(X,c_{X,\_})$, where $X$ is an object of $\C$ and $c_{X,\_}$ is a family of natural isomorphisms $c_{X,Y} : X \otimes Y \rightarrow Y \otimes X$ such that, for all $Y$ and $Z$, the following diagram commutes:
 $$
\xymatrix{X \otimes Y \otimes Z \ar[rr]^{c_{X,Y \otimes Z}} \ar[dr]^{c_{X,Y} \otimes Z} & & Y \otimes Z \otimes X\\
 & Y \otimes X \otimes Z \ar[ur]^{Y \otimes c_{X,Z}} &
 }
 $$
If $\C$ is a fusion category, $\mathcal{Z}(\C)$ is a modular category, and in particular a braided fusion category. Moreover, the Frobenius-Schur indicators of simple objects of $\mathcal{Z}(\C)$ can be extracted from its S and T matrices.

A \emph{pointed fusion category}  is a fusion category  in which every object is invertible with respect to the tensor product. Pointed fusion categories are parametrized by couples $(G,\omega)$ where $G$ is a finite group and $\omega$ is a $\mathbb{C}^\times$-valued 3-cocycle on $G$ - they are denoted ${\mathcal{C}}(G,\omega)$. Simple objects of these are the elements of $G$, the tensor product is given by multiplication and the associativity constraint between $g,h,k \in G$ is the scalar tranformation given by $\omega(g,h,k)\Id : (gh)k \to g(hk)$. We note that the category $\operatorname{Vect}^\omega_G$ of $G$-graded vector space with twisted tensor product $(u \otimes v) \otimes w)=\omega(g,h,k) u \otimes (v \otimes w)$ for respectively $g$, $h$ and $k$ homogenous elements $u$,$v$ and $w$. is often used.  $\C(G,\omega)$ is a skeletal category equivalent $\operatorname{Vect}^\omega_G$.  We refer the reader to \cite{O}.

Following \cite{Ni}, we say that two fusion categories ${\mathcal{C}}$ and ${\mathcal{D}}$ are called categorically Morita equivalent if and only if there exists an indecoposable ${\mathcal{C}}$-module category ${\mathcal{M}}$ such that ${\mathcal{D}}$ is equivalent (as fusion categories) to ${\mathcal{C}}_{\mathcal{M}} ^*$. This is an equivalence relation that contains the usual equivalence of categories. Also, two categories are categorically Morita equivalent if and only if they have equivalent (as braided fusion categories) centers. Categories that are Morita equivalent to pointed fusion categories are called \emph{group theoretical categories} (g.t.c.). They are parametrized by quadruples $(G,\omega,H,\mu)$, where $G$ is a finite group, $H$ is a subgroup of $G$, $\omega$ a $3$-cocycle on $G$ and $\mu$ a 2-cochain on $H$ such that $d \mu = \omega|_{H \times H \times H}$. Such a quadruple is called a \emph{group theoretical data} and the associated category is denoted ${\mathcal{C}}(G,\omega,H,\mu)$.

\subsection{Projective characters}
\label{sec:projc}

Simple objects of group-theoretical categories are parametrized by irreducible projective characters of certain subgroups of the groups defining the categories. Since we need to do explicit calculations with these objects, we recall here some basic notions from projective character theory, and we give the GAP codes we use for dealing with them.  General references for projective representation theory are \cite{Karpilovsky} or  \cite{Isaacs}.

A projective representation is a group homomorphism $\underline\rho\colon G\to\PGL(V)$ for a finite dimensional vector space $V$. If we choose an inverse image $\rho(g)\in\GL(V)$ for $\underline\rho(g)\in\PGL(V)=\GL(V)/Z(\GL(V))$ for each $g\in G$ (taking $\rho(1)=1$), then we obtain a map $\rho\colon G\to\GL(V)$ which is a group homomorphism only up to scalar factors:
\begin{definition}
  Let $\alpha: G \times G \rightarrow \mathbb{C}^\times$ be a map. An \emph{$\alpha$-projective representation}, or simply \emph{$\alpha$-representation}, of $G$, is the data of a complex vector space $V$ and a mapping $\rho : G \rightarrow GL(V)$, such that the following equalities hold:
\begin{eqnarray}
\rho(gh)& = & \rho(g) \rho(h) \alpha^{-1}(g,h)\label{2}\\
\rho(1) & = & 1
\end{eqnarray}
The character of $\rho$ is the map $\tr\circ \rho\colon G\to\CC$, and an $\alpha$-projective character is the character of an $\alpha$-projective representation. 
\end{definition}
If $\rho$ is an $\alpha$-projective representation, then $\alpha$ is necessarily a normalized two-cocycle.

Consider the central extension 
$$1 \rightarrow \mu_n \rightarrow G_{\alpha} \rightarrow G \rightarrow 1$$
defined by such a two-cocycle, where $G_{\alpha}$ denotes the ``twisted direct product'' $\mu_n \times_{\alpha} G$ with the group $\mu_n$ of $n$-th roots of unity, whose structure is given by:
$$(z , g) \cdot (z',g')=(z z' \alpha(g,g'),g g') \ , \ \forall z,z' \in \mu_n \ , \ \forall g,g' \in G$$

Then any $\alpha$-projective representation $\rho$ gives rise to an ordinary representation $\tilde\rho\colon G_\alpha\to\GL(V)$ by $\tilde\rho(z,g)=z\rho(g)$, whose restriction to $\mu_n$ is the tautological representation of $\mu_n$ by multiplication. Conversely, any representation of $G_\alpha$ whose restriction to $\mu_n$ is the tautological representation of $\mu_n$ arises in this fashion from an $\alpha$-representation of $G$. Note that, for a representation $\rho_\alpha$ of $G_\alpha$ and its character $\chi_\alpha$, the restriction of $\rho_\alpha$ to $\mu_n$ is the tautological representation iff the restriction of $\chi_\alpha$ to $\mu_n$ is the character of the tautological representation, that is $\chi_\alpha(\zeta_n)=\dim(V)\zeta_n$ for some generator $\zeta_n$ of $\mu_n$.

We now give our implementation in GAP of $\alpha$-representations through central extensions.
We use the package HapCocyclic of GAP for this implementation. We did add a few new commands and modified some code to improve the efficiency of the package for our purposes, but we will not discuss these merely technical details here.
\lstinputlisting[linerange=proj-end]{global.txt}\label{lst:projcode}
This will allows us to compute with GAP simple objects of a group theoretical category ${\mathcal{C}}(G,\omega,H,1)$ with adapted $\omega$, that are parametrized by representatives $g$ of right cosets $gH$ and $\omega_g$-projective characters of the stabilizer $S:=Stab_H(gH)=H \bigcap g \triangleright H$.

\subsection{Frobenius-Schur indicators}
\label{ssec:FSI}
It is easy to modify general group theoretical data $(G,\omega,H,\mu)$ to give data of the form $(G,\omega',H,1)$ that gives rise to an equivalent fusion category. One simply has to divide $\omega$ by the coboundary of an arbitrary extension of $\mu$ to a cochain on $G$. In \cite{Na}, Natale showed that one can even replace $\omega$ by $\omega'$ in such a way that $\omega'|_{G \times G \times H}$ is trivial, which is a slightly stronger condition. We call such $\omega$ \emph{adapted} cocycle. Using that trick, one of the authors \cite{S1,Serratum} extracted formulas for higher indicators of simple objects in group-theoretical categories, for both adapted and general cocycles. As a special case, a formula for doubles is also given. We recall those formulas and refer to \cite{S1,Serratum} for further details. The way we compute adapted cocycles is the same as in \cite{Serratum} and we give its implementation at the end of this section.

We define for $g,x,y \in G$ and $\omega \in Z^3(G,\mathbb{C}^\times)$ \ the following symbols:
\begin{align*}
\omega_g(x,y)&:=\omega(x,y,g)\\
\alpha_g(x,y)&:=\omega(x,y,g)\omega^{-1}(x,y \triangleright g,y)\omega(xy \triangleright g,x,y),
\end{align*}
where $\triangleright$ denotes the conjugation action.
We also define the symbols $\pi_m(x)$ fixed for $m\in\ZZ$ by
$$\pi_0(x):=1 \text{ and } \pi_{m+1}(x):=\omega(x,x^m,x)\pi_m(x)\text{ for }m\in\ZZ$$
As the sequence of numbers $\omega(x,x^m,x)$ is periodic of period $o$ the order of $x \in G$, we can use the following to compute efficiently $\pi_m(x)$:
\begin{align}
\pi_{m+o}(x) 
& = \left(\prod_{m \leqslant k \leqslant m+o-1} \omega(x,x^k,x)\right) \pi_m(x)\nonumber\\
& = \left(\prod_{0 \leqslant k \leqslant o-1} \omega(x,x^k,x)\right) \pi_m(x)\nonumber\\
& =  \pi_o(x) \pi_m(x)\label{eq:pi}
\end{align}

\begin{theorem}\label{GTformula}
For a simple object $M$ in ${\mathcal{C}}(G,\omega,H,1)$ with adapted cocycle corresponding to the couple $(g,\chi)$ where $\chi$ is an $\omega_g$-projective character of $S$, we have:
\begin{eqnarray}
\nu_m(M) &=& \frac{1}{|S|}\sum_{
\begin{tabular}{c}
$x \in gH$\\
 $x^m \in S$
\end{tabular}}
\pi_{-m}(x)\chi(x^{-m})\\
 &=& \frac{1}{|S|}\sum_{
\begin{tabular}{c}
$h \in H$\\
 $(gh)^m \in S$
\end{tabular}}
\pi_{-m}(gh)\chi((gh)^{-m})
\end{eqnarray}
\end{theorem}
\begin{cor}\label{PTformula}
For a simple object $M$ in ${\mathcal{C}}(G,\omega)$ corresponding to $g$, we have:
$$\nu_m(M) =
\left\{
\begin{tabular}{cc}
$ \pi_{-m}(g)$& if\ \ $g^m=1_G$\\
$ 0$ &  if\ \ $g^m \neq 1_G$
\end{tabular}
\right.
$$
\end{cor}
\ 

Now, the simple objects of the category $\mathcal{Z}(\mathcal{C}(G,\omega))$ are parametrized by representatives $g$ of conjugacy classes of $G$ and $\alpha_g$-projective characters on $C_G(g)$.
\begin{theorem}\label{doublesformula}
For a simple object $M$ in $\mathcal{Z}(\mathcal{C}(G,\omega))$ corresponding to the couple $(g,\chi)$ where $\chi$ is an $\alpha_g$-projective character of $C_G(g)$, we have:
\begin{equation}
\nu_m(M) = \frac{1}{|C_G(g)|}\sum_{
\begin{tabular}{c}
$x \in G$\\
 $(gx)^m = x^m$
\end{tabular}}
\alpha_{x^m}(g,x)\frac{\pi_m(gx)}{\pi_m(x)}\chi(x^m)
\end{equation}
\end{theorem}

The higher indicators are invariants of pivotal fusion categories. In particular, they are not related to the braiding in $\mathcal{Z}(\mathcal{C}(G,\omega))$. Some information about this aspect of the double is contained in the ribbon structure. 

Recall from \cite{GMN} the following:
\begin{lemma}\label{doublestwist}
For a simple object $M$ in $\mathcal{Z}(\mathcal{C}(G,\omega))$ associated to the couple $(g,\chi)$, the ribbon structure $\theta_M \in Hom(M,M^*)\cong \mathbb{C}^\times$ is given by the formula:
$$\theta_M=\chi(g)/\chi(1_G)$$
\end{lemma}
Finally, we recall that the sequence of Frobenius-Schur indicators of a simple object in a fusion category is periodic. The least common multiple of those periods for all simple object is called the \emph{Frobenius-Schur exponent} of the category and is studied in \cite{NS2}. The authors showed in particular that this exponent is exactly the order of the ribbon twist of the center. For a pointed category $\mathcal{C}(G,\omega)$, the FS exponent is given by:
$$\FSexp( \mathcal{C}(G,\omega))=\lcm(|C||\omega_C|)$$
where $C$ runs through the cyclic subgroups of $G$ and $|\omega_C|$ is the order of the cohomology class of the restriction of $\omega$ to $C$.

We finish this section by our implementation in GAP of the formulas given in \cref{PTformula,GTformula,doublesformula,doublestwist}.

We use functoriality of cohomology groups to get the exponent of the category $\C(G,\omega)$, for a group $G$, a free resolution $R$ for $G$, the complex $K:=R \otimes \mathbb{Z}$, the record given by our function \lstinline!UniversalCoefficientsTheorem(K,3)!  and a $3$-cocycle $\omega$:
\lstinputlisting[linerange=exponent-end]{global.txt}\label{lst:expcode}

Now, we will use cocycles in the bar resolution. For an \lstinline!n!-cocycle \lstinline!omega! in a resolution \lstinline!R! with values in a cyclic ring of order \lstinline!m!, we use
\begin{lstlisting}
StandardCocycle(R,omega,n,m);
\end{lstlisting}
for getting the image of \lstinline!omega! along a chain map between \lstinline!R! and the bar resolution. we note that for using this command, HAP needs a contracting homotopy of \lstinline!R!.

The symbols $\pi_m$ we defined in \ref{ssec:FSI} can be computed, for an element $x \in G$ of order $o$, $m \in \{0, \dots, o\}$ and a $3$-cocycle $\omega \in H^3(G,\mathbb{C}^\times)$ of exponent $n$, by the following function. In order to get $\pi_m(x)$ for values of $m$ that are negative or greater than $o$, we use \cref{eq:pi}. For convenience, the code stores the values for $m \in \{0 \dots e\}$, where $e:=\Exp(G)$, for every element of $G$.
\lstinputlisting[linerange=pi-end]{global.txt}\label{lst:picode}

We can then give the functions that compute our invariants. Again, we describe simple objects of $\mathcal{Z}(\C(G,\omega))$ by couples $(g,\chi)$ where $g$ is a representative of a conjugacy class of $G$, and $\chi$ is an $\alpha_g$-projective character of $C_G(g)$. 

First, the function \lstinline!Twist! gives the ribbon structure for a simple object of $D_\omega(G)$. The arguments of this function are: \lstinline!class! represents the representative $g$, \lstinline!C! is the centralizer $C_G(g)$, lift is the section of the epimorhism $C_{\alpha_g} \rightarrow C$ and finally \lstinline!chi! is a character of $C_{\alpha}$ induced by an $\alpha_g$-projective character of \lstinline!C!.
\lstinputlisting[linerange=twist-end]{global.txt}\label{lst:twistcode}

Now, we store all the values of the $\pi$ symbols for all $x \in G$ in a list.
\begin{lstlisting}
listG:=Enumerator(G);
pivalues:=List(listG,x->
	ValuesOfPiSymbols(G,x,n,f));
\end{lstlisting}

Finally, the following function gives the $m$-Frobenius-Schur indicator of a simple object.

\lstinputlisting[linerange=indicator-end]{global.txt}\label{lst:doublesFSIcode}

And also, we give the codes that computes the indicators of pointed fusion categories and adapted group-theoretical categories. First, for pointed categories.
\lstinputlisting[linerange=ptindicators-end]{global.txt}\label{lst:ptcode}
Now we show how, from a group-theoretical data $(G,\omega,H,1)$, we adapt the $3$-cocycle $\omega$ whose restriction to $H$ is trivial to a $3$-cocycle $\omega'$ such that $\omega'|_{G \times G \times H}=1$.
\lstinputlisting[linerange=adaptation-end]{global.txt}\label{lst:adaptation}
And finally, we can give the code for the FSI of group-theoretical categories with adapted cocycles.
\lstinputlisting[linerange=gtindicators-end]{global.txt}\label{lst:gtcode}
\section{Classification of $\mathcal{C}(G,\omega)$ up to $|G|<32$}
\label{sec:work}
In \cite{GMN}, the authors showed that there exist exactly 20 Morita equivalence classes of $\mathcal{C}(G,\omega)$ for $\vert G \vert=2^3=8$, but only for the sets of couple $(G,\omega)$ for whom the bialgebra $D_{\omega}(G)$ is non commutative. In order to obtain this result, they went through the following strategy. First, find representatives for the orbits of $3$-cocycles under automorphisms of $G$. Then, find Morita equivalences coming from \cref{theo21}. This gives an upper bound for the number of classes which is lower that just the number of automorphism orbits of cocycles. Finally compute the FS indicators and twists associated to chosen representatives of couples $(G,\omega)$ under the equivalence relation generated by \cref{theo21} and automorphism orbits. Since the number of different sets of table is the same as the number of those representatives, they lead to non-equivalent categories. 

We will here follow a slightly different strategy, which uses computer help more directly and needs less knowledge about the cohomology groups considered. We will not distinguish the cases where the bialgebra structure of $D_{\omega}(G)$ is commutative or not, but take for a fixed order of $G$ all pointed categories $\mathcal{C}(G,\omega)$. We take automorphism group orbits of cocycles since we know that cocycles in the same orbit give equivalent categories. Then we collect all those categories that give the same tables of invariants, and see if we can show that the categories affording the same invariants are in fact equivalent.

First, we describe our approach for some particular orders, mostly for $|G|=p^\alpha$ where p is a prime. This includes the known case where $G$ has order $2^3=8$, where \ref{corcentral}, together with computer calculations, is sufficient to conclude; we obtain then similar results for the case $|G|=3^3=27$. Finally, we treat the more complicated case $|G|=2^4=16$, for which we have to use indicators of some group-theoretical categories to conclude. We also not that Uribe et Munoz \cite{Uribe} obtained the same results for $|G|=8$ using the Lyndon-Hochschild-Serre spectral sequence of finite groups. 
\subsection{Morita equivalence of pointed fusion categories}
\label{ssec:Morita}

Equivalences of two pointed fusion categories are described as follows: Two pointed categories $(G_1,\omega_1)$ and $(G_2,\omega_2)$ are equivalent if and only if there exists an isomorphism $\phi:G_1 \cong G_2$ such that $\omega_1$ is cohomologous to $\phi^\ast (\omega_2)$. In particular, consider the action of the group $\Phi$ of automorphisms of $G$ on $H^3(G,\CCu)$ given by $(\omega \triangleleft \phi)(g,h,k):=\omega(\phi(g),\phi(h),\phi(k))$ with $\omega \in Z^3(G,\mathbb{C}^\times)$ and $\phi \in \Phi$. Then two categories $\C(G,\omega)$ and $\C(G,\omega')$ are equivalent iff the cohomology classes of $\omega$ and $\omega'$ are in the same automorphism orbit. 

Trivially, equivalent categories are categorically Morita equivalent. In particular, $\mathcal{C}(G_1,\omega_1)$ and $\mathcal{C}(G_2,\omega_2)$ are categorically Morita equivalent if there is an isomorphism $\phi\colon G_1\to G_2$ with $\omega_1$ and $\phi^*\omega_2$ cohomologous, and if $G_2=G_1=G$, they are equivalent if the cohomology classes of $\omega_1$ and $\omega_2$ are in the same orbit under automorphisms of $G$. 

For a group $G$ and its automorphisms group $\Phi$, we give a procedure to obtain automorphism orbits of $H^n(G,\mathbb{C}^\times)$ under the action of $\Phi$. This gives a list of lists, whose elements are the different automorphism orbits. For this, we use a part of the functoriality described in \ref{ssec:GAP1}.
\lstinputlisting[linerange=automorphisms-end]{global.txt}\label{lst:autocode}

To treat categorical Morita equivalence of pointed fusion categories, which is a larger equivalence relation, one needs to understand when $\C(G,\omega,H,\mu)$ is pointed, and identify the group $G'$ and the cocycle $\omega'$ in a category equivalence $\C(G,\omega,H,\mu)\cong\C(G',\omega')$. In principle, an approach to solving this is contained in \cite{S2}, where for an adapted cocycle, a coquasibialgebra whose comodule category is $\C(G,\omega,H,1)$ is explicitly described; however, the description is not as computationally accessible as to allow an immediate extraction of $G'$ and $\omega'$ after deciding if the category is pointed (i.e. the coquasibialgebra cocommutative). In a special case that we will discuss below, an explicit description of $(G',\omega')$ is in fact extracted in \cite{GMN}. As we will point out, it is particularly convenient since, in conjunction with our use of invariants of the centers, it allows to deduce certain categorical Morita equivalences by calculations on the level of cohomology classes (rather than standard cocycles).  In the general situation, a complete solution is in \cite{Deepak} in terms of group cohomology. This solution involves calculations with cocycles in the standard resolution with nontrivial coefficients, however a more cohomological interpretation of the results of \cite{Deepak} is given in \cite{Uribe1} and applied to the case of groups of order eight in \cite{Uribe}.

We did not attempt to treat the criteria from  \cite{Deepak} or \cite{Uribe1} computationally, but used more ad hoc methods to find categorical Morita equivalences. The first one is a special case of Theorem \ref{theo21} of \cite{GMN} stated in Corollary \ref{corcentral}; the second one involves computation of Frobenius-Schur indicators and is described in our examples, see \ref{sssec:16}.

{We recall now the result of \cite{GMN}}. Take $H$ a finite group and $A$ a right $H$-module with action $\triangleleft$. We define two different extensions of $H$. The first one is the semi-direct product $E:=H \ltimes A$ given by multiplication 
$$(x,a)(y,b):=(xy,(a \triangleleft y)b)$$
Then we have the following extension
$$1 \rightarrow A \rightarrow E \rightarrow H \rightarrow 1$$
For the second extension, take $\hat{A}=Hom(A,\mathbb{C}^\times)$ with left $H$-action $\triangleright$:
$$(h \triangleright \chi)(a):=\chi(a \triangleleft h)$$
Then, for any 2-cocycle $\eta \in Z^2(H,\hat{A})$, one can construct an extension 
$$1 \rightarrow \hat{A} \rightarrow G \rightarrow H \rightarrow 1$$
where $G:=\hat{A} \rtimes_ {\eta} H$ has multiplication
$$(\chi ,x)(\psi ,y):=(\chi(x \triangleright \psi)\eta(x,y),xy)$$
For $\zeta \in Z^3(H,\mathbb{C}^\times)$, using the projections of both extensions, we can define the \emph{inflations} $\zeta_E \in Z^3(E,\mathbb{C}^\times)$ and $\zeta_G \in Z^3(G,\mathbb{C}^\times)$ of $\zeta$. Also, for $\sigma$ a $H$-module automorphism of $A$, one can define the 3-cocycle $\omega_{\sigma} \in Z^3(E,\mathbb{C}^\times)$ by
$$\omega_ {\sigma} \left( (a_1,h_1),(a_2,h_2),(a_3,h_3)\right):= \eta(h_2,h_3) \left(\sigma(a_1) \right)$$

\begin{theorem}\label{theo21}
Take $\sigma \in Aut_H (A)$. The categories $\mathcal{C}(G,\zeta_G)$ and $\mathcal{C}(E,\omega_\sigma \zeta_E)$ are categorically Morita equivalent.
\end{theorem}

  \begin{cor}
    \label{corcentral} Let $1\to A\to G\to Q\to 1$ be an extension with $A$ abelian, and $\eta\in\in H^3(Q,\CCu)$. Then there is a semidirect product $A\semdir Q$ and a three-cocycle $\alpha$ on $A\semdir Q$ such that $\C(G,\Inf_Q^G\eta)$ is categorically Morita equivalent to $\C(A\semdir Q,\alpha)$. If the extension $G$ is central, the semidirect product is a direct product.
  \end{cor}
\subsection{Some cases with $|G|=p^\alpha$}
\label{ssec:prime}
\subsubsection{$|G|=2^3$}
\label{sssec:8}
First, we review the basic facts on finite groups of order $8$. The structure of the group $G$, the cohomology group $H^3(G,\mathbb{C})$, the size of the group $Aut(G)$ of automorphism of $G$ and the number of automorphisms orbits of $3$-cocycles are summarized in the following table:
\begin{center}
\begin{tabular}{|c|c|c|c|}

\hline
$G$ & $H^3(G,\mathbb{C}^\times)$ & $|Aut(G)|$ & Number of orbits\\
\hline
$\mathbb{Z}/8 \mathbb{Z}$ & $\mathbb{Z}/8 \mathbb{Z}$ & $4$ & $8$\\
$\mathbb{Z}/4 \mathbb{Z} \times \mathbb{Z}/2 \mathbb{Z}$ & $\mathbb{Z}/4 \mathbb{Z} \times \mathbb{Z}/2 \mathbb{Z} \times \mathbb{Z}/2 \mathbb{Z}$ & $8$ & $9$\\
$D_8$ & $\mathbb{Z}/4 \mathbb{Z} \times \mathbb{Z}/2 \mathbb{Z} \times \mathbb{Z}/2 \mathbb{Z}$ & $8$ & $12$\\
$Q_8$ & $\mathbb{Z}/8 \mathbb{Z}$ & $24$ & $8$\\
$(\mathbb{Z}/2 \mathbb{Z})^3 $ & $(\mathbb{Z}/2 \mathbb{Z})^7$ & $168$ & $10$\\
\hline
\end{tabular}
\captionof{table}{}\label{tab:groups8}
\end{center}

Thus there are at most $47$ Morita equivalence classes of pointed fusion categories of rank $8$. We obtained a lower bound of the number of those classes by computing our invariants. Let us first fix some notations. The five groups of order $8$ are labeled by integers from $1$ to $5$ in the same order as in the table; for example, the dihedral group is denoted as $G_3$. For each group, we fix an order on the set of representatives of automorphism orbits of cocycles we are given; the $j$-th representatives on the $i$-th group is denoted as $\omega_i^j$.  For the $47$ objects we consider, we obtain $38$ different sets of invariants, which gives us a lower bound. Now, we prove that this bound is optimal: We consider the different pairs of a group and a cocycle that share the same tables of invariants, to see whether the associated twisted doubles are in fact equivalent.

Among the list of $38$ tables, there are already $30$ that appear only once, so we know that those characterize their respective gauge classes. For the eight remaining tables, seven appear twice and one appears three times. For example, the following table comes from the couples $(G_2, \omega_2^7)$ and $(G_5 ,\omega_5^5)$: The first columns labelled $\mu_i$ contain the values of FS indicators, the column labelled $\theta$ gives the ribbon twists and the last colum gives the number of simple objects that share a same set of those invariants.

\begin{center}
\begin{tabular}{|c|c|c|c|c|}
\hline
$ \nu_0$ & $\nu_1$ & $\nu_2$ & $\theta$ & $m$\\
\hline
$1$ & $1$ &$1$ & $1$ &$1$\\
$1$ & $0$ &$0$ & $1$ &$12$\\ 
$1$ & $0$ &$1$ & $1$ &$7$\\
$1$ & $1$ &$1$ & $-1$ &$1$\\
$1$ & $0$ &$0$ & $i$ &$12$\\
$1$ & $0$ &$0$ & $-i$ &$12$\\
$1$ & $0$ &$-1$ & $i$ &$12$\\
$1$ & $0$ &$-1$ & $-i$ &$12$\\
\hline
\end{tabular}
\captionof{table}{FS indicators for $\mathcal{Z}(\mathcal{C}(G_2,\omega_2^7))$ and $\mathcal{Z}(\mathcal{C}(G_5,\omega_5^5))$}\label{tab:FSIcouple8}
\end{center}
The table that is repeated three times is the one for the couples $(G_3,\omega_3^5)$, $(G_4,\omega_4^1)$ and $(G_5, \omega_5^4)$:
\begin{center}
\begin{tabular}{|c|c|c|c|c|}
\hline
$ \nu_0$ & $\nu_1$ & $\nu_2$ & $\theta$ & $m$\\
\hline
$1$ & $1$ &$1$ & $1$ &$1$\\
$1$ & $0$ &$1$ & $1$ &$7$\\
$2$ & $0$ &$1$ & $1$ &$3$\\
$2$ & $0$ &$1$ & $-1$ &$3$\\
$2$ & $0$ &$-1$ & $1$ &$1$\\
$2$ & $0$ &$-1$ & $-1$ &$1$\\
$2$ & $0$ &$-1$ & $i$ &$3$\\
$2$ & $0$ &$-1$ & $-i$ &$3$\\
\hline
\end{tabular}
\captionof{table}{FS indicators for $\mathcal{Z}(\mathcal{C}(G_3,\omega_3^5))$, $\mathcal{Z}(\mathcal{C}(G_4,\omega_4^1))$ and $\mathcal{Z}(\mathcal{C}(G_5,\omega_5^4))$}\label{tab:FSItriple8}
\end{center}

Now, we explain how to deal with the seven occurrences of a table repeated twice. We list the couples $(G_i,\omega_i^j)$ involved by pairing them according to their tables:\\
$$(G_2,\omega_2^7) \text{ and } (G_5, \omega_5 ^5)\text{ ,}$$
$$(G_2,\omega_2^1) \text{ and } (G_5, \omega_5 ^3)\text{ ,}$$
$$(G_1,\omega_1^1) \text{ and } (G_2, \omega_2 ^2)\text{ ,}$$
$$(G_1,\omega_1^5) \text{ and } (G_2, \omega_2 ^4)\text{ ,}$$
$$(G_4,\omega_4^2) \text{ and } (G_5, \omega_5 ^6)\text{ ,}$$
$$(G_3,\omega_3^9) \text{ and } (G_5, \omega_5 ^8)\text{ ,}$$
$$(G_3,\omega_3^1) \text{ and } (G_5, \omega_5 ^2)\text{ .}$$

We discuss the first case, the other ones can be treated in the same way. We use Corollary \ref{corcentral}. For the groups $G_2=\mathbb{Z}/ 4\mathbb{Z} \times \mathbb{Z}/2 \mathbb{Z}$ and $G_5=(\mathbb{Z}/2\mathbb{Z})^3$, there is a common abelian normal subgroup $A$, isomorphic to $\mathbb{Z}/2\mathbb{Z}$, such that the inflation of a $3$-cocycle on $Q \cong (\mathbb{Z}/2 \mathbb{Z})^2$ to $G_2$ lies in the same cohomology class as $\omega_2^7$. We find this cocycle by using the functoriality along the natural projection of $G_5$ over $Q$, and computing the inflations in $H^3(G,\mathbb{C}^\times)$ of all cohomology classes of $H^3(Q,\mathbb{C}^\times)$. The theorem therefore asserts that there must be a $3$-cocycle $\alpha$ on $G_5$ such that $\mathcal{C}(G_2,\omega_2^7)$ and $\mathcal{C}(G_5,\alpha)$ are Morita equivalent. We can conclude that $\alpha$ must be in the orbit of  $\omega_5^5$, since we have calculated invariants for each cocycle orbit on $G_5$, and  only the invariants of $\mathcal{Z}(\mathcal{C}(G_5,\omega_5^5))$ agree with those of $\mathcal{Z}(\mathcal{C}(G_2,\omega_2^7))$, as they have to if the categories are to be equivalent.

We still have to deal with the table that appears for the couples $(G_3,\omega_3^5)$, $(G_4,\omega_4^1)$ and $(G_5, \omega_5^4)$. Despite the fact that there are three couples that share this table, the same argument as before works  since $G_5$ has only one cocycle orbit affording the same invariants as the other two. Thus, we can find Morita equivalences between $\mathcal{C}(G_3,\omega_3^5)$ and $\mathcal{C}(G_5,\omega_5^4)$ first, and then between $\mathcal{C}(G_4,\omega_4^1)$ and $\mathcal{C}(G_4,\omega_5^4)$ (which is enough since Morita equivalence is an equivalence relation).

More precisely, both of the groups $D_8$ and $Q_8$ have a central subgroup isomorphic to $\mathbb{Z}/2\mathbb{Z}$, with the quotient isomorphic to $(\mathbb{Z}/2\mathbb{Z})^2$ (and so, the direct product is $(\mathbb{Z}/2\mathbb{Z})^3$). It is easy to deduce a Morita equivalence between $\mathcal{C}(G_4,\omega_4^1)$ and $\mathcal{C}(G_5,\omega_5^4)$, as the cocycle $\omega_4^1$ is the trivial $3$-cocycle on $Q_8$ and is therefore the inflation of the trivial $3$-cocycle on the quotient. For the equivalence between $\mathcal{C}(G_3,\omega_3^5)$ and $\mathcal{C}(G_5,\omega_5^4)$, we find again a $3$-cocycle on $(\mathbb{Z}/2\mathbb{Z})^2$ with inflation $\omega_3 ^5$. We can then state the following:
\begin{theorem}\label{theo8}
There are exactly $38$ Morita equivalence classes of pointed fusion categories of rank $8$.
\end{theorem}
\subsubsection{$|G|=3^3$}
\label{sssec:27}
As in the previous subsection, we recall all information we need about groups of order $27$, their cohomology, and the number of orbits under the action of the automorphism group.
$$\begin{tabular}{|c|c|c|c|}
\hline
$G$ & $H^3(G,\mathbb{C}^\times)$ & $|Aut(G)|$ & Number of orbits\\
\hline
$\mathbb{Z}/27 \mathbb{Z}$ & $\mathbb{Z}/27 \mathbb{Z}$ & $18$ & $7$\\
$\mathbb{Z}/9 \mathbb{Z} \times \mathbb{Z}/3 \mathbb{Z}$ & $\mathbb{Z}/9 \mathbb{Z} \times \mathbb{Z}/3 \mathbb{Z} \times \mathbb{Z}/3 \mathbb{Z}$ & $108$ & $16$\\
$(\mathbb{Z}/3 \mathbb{Z})^2 \rtimes \mathbb{Z}/3 \mathbb{Z}$ & $(\mathbb{Z}/3 \mathbb{Z})^4$ & $432$ & $15$\\
$\mathbb{Z}/9 \mathbb{Z} \rtimes \mathbb{Z}/3 \mathbb{Z}$ & $(\mathbb{Z}/3 \mathbb{Z})^2$ & $54$ & $9$\\
$(\mathbb{Z}/3 \mathbb{Z})^3 $ & $(\mathbb{Z}/3 \mathbb{Z})^7$ & $11232$ & $14$\\
\hline
\end{tabular}$$
We then compute the $61$ tables of our invariants for the given representatives of automorphism orbits. There are $47$ different tables, so some appear for several group-cocycle pairs: $34$ tables appear once, $12$ tables appear twice, and one table appears three times. The couples involved are given, with the same kind of notations as before:
$$(G_1,\omega_1^7) \text{ and } (G_2, \omega_2 ^6)\text{ ,}$$
$$(G_1,\omega_1^6) \text{ and } (G_2, \omega_2 ^4)\text{ ,}$$
$$(G_1,\omega_1^1) \text{ and } (G_2, \omega_2 ^2)\text{ ,}$$
$$(G_2,\omega_2^9) \text{ and } (G_5, \omega_5 ^5)\text{ ,}$$
$$(G_2,\omega_2^{14}) \text{ and } (G_5, \omega_5 ^7)\text{ ,}$$
$$(G_2,\omega_2^1) \text{ and } (G_5, \omega_5 ^3)\text{ ,}$$
$$(G_3,\omega_3^1) \text{ and } (G_5, \omega_5 ^2)\text{ ,}$$
$$(G_4,\omega_4^4) \text{ and } (G_5, \omega_5 ^6)\text{ ,}$$
$$(G_4,\omega_4^7) \text{ and } (G_5, \omega_5 ^8)\text{ ,}$$
$$(G_3,\omega_3^{13}) \text{ and } (G_5, \omega_5 ^{14})\text{ ,}$$
$$(G_3,\omega_3^{10}) \text{ and } (G_5, \omega_5 ^{12})\text{ ,}$$
$$(G_3,\omega_3^7) \text{ and } (G_5, \omega_5 ^{10})\text{ ,}$$
$$(G_3,\omega_3^5) \text{ , } (G_4,\omega_4^1)\text{ and } (G_5, \omega_5 ^{4})\text{.}$$
In the same way as in section \ref{sssec:8}, using Corollary \ref{corcentral} , we find Morita equivalences between all the $12$ pairs that share a table, and also between $\mathcal{C}(G_3,\omega_3^5)$ and $\mathcal{C}(G_5,\omega_5^4)$ and $\mathcal{C}(G_4,\omega_4^1)$ and $\mathcal{C}(G_5,\omega_5^4)$. And so, we can state:
\begin{theorem}\label{theo27}
There are exactly $47$ Morita equivalence classes of pointed fusion categories of rank $27$.
\end{theorem}
\subsubsection{$|G|=2^4$}
\label{sssec:16}
There are $14$ groups of order $16$ for $17.312$ pairs of a group and a cohomology class. As usual, we gathered all this data in the next table, along with the orders of the automorphisms group and the number of different classes of $3$-cocycles under the action of the automorphism group. One can note that the third cohomology of the abelian group $(\mathbb{Z}/2 \mathbb{Z})^4$ is quite large (it is a vector space of dimension $14$ over $\mathbb{Z}/2\mathbb{Z}$), but luckily its automorphism group is also large ($\operatorname{GL}_4(2)$, of order 20.160). In the end, there are only $23$ automorphism orbits to consider for this group.\\
\begin{center}
\begin{tabular}{|c|c|c|c|}
\hline
$G$ & $H^3(G,\mathbb{C}^\times)$ & $|Aut(G)|$ & \begin{tabular}{c}Number \\of orbits\end{tabular}\\
\hline
$\mathbb{Z}/16 \mathbb{Z}$ & $\mathbb{Z}/16 \mathbb{Z}$ & $8$ & $12$\\
$\mathbb{Z}/4 \mathbb{Z} \times \mathbb{Z}/4 \mathbb{Z}$ & $(\mathbb{Z}/4 \mathbb{Z})^2$ & $96$ & $12$\\
$(\mathbb{Z}/4 \mathbb{Z} \times \mathbb{Z}/2 \mathbb{Z}) \rtimes \mathbb{Z}/2 \mathbb{Z} $ & $(\mathbb{Z}/2 \mathbb{Z})^2 \times (\mathbb{Z}/4 \mathbb{Z})^2$ & $32$ & $30$\\
$\mathbb{Z}/4 \mathbb{Z} \rtimes \mathbb{Z}/4 \mathbb{Z}$ & $\mathbb{Z}/2 \mathbb{Z} \times (\mathbb{Z}/4 \mathbb{Z})^2$ & $32$ & $24$\\
$\mathbb{Z}/8 \mathbb{Z} \times \mathbb{Z}/2 \mathbb{Z}$ & $(\mathbb{Z}/2 \mathbb{Z})^2 \times \mathbb{Z}/8 \mathbb{Z}$ & $16$ & $18$\\
$\mathbb{Z}/8 \mathbb{Z} \rtimes \mathbb{Z}/2 \mathbb{Z}$ & $\mathbb{Z}/2 \mathbb{Z} \times \mathbb{Z}/8 \mathbb{Z}$ & $16$ & $16$\\
$D_{16}$ & $(\mathbb{Z}/2 \mathbb{Z})^2 \times \mathbb{Z}/8 \mathbb{Z}$ & $32$ & $24$\\
$QD_{16}$ & $\mathbb{Z}/2 \mathbb{Z} \times \mathbb{Z}/8 \mathbb{Z}$ & $16$ & $16$\\
$Q_{16}$ & $\mathbb{Z}/8 \mathbb{Z}$ & $32$ & $12$\\
$\mathbb{Z}/4 \mathbb{Z} \times (\mathbb{Z}/2 \mathbb{Z})^2$ & $(\mathbb{Z}/2 \mathbb{Z})^6 \times \mathbb{Z}/4 \mathbb{Z}$ & $192$ & $34$\\
$\mathbb{Z}/2 \mathbb{Z} \times D_8$ & $(\mathbb{Z}/2 \mathbb{Z})^6 \times \mathbb{Z}/4 \mathbb{Z}$ & $64$ & $57$\\
$\mathbb{Z}/2 \mathbb{Z} \times Q_8$ & $(\mathbb{Z}/2 \mathbb{Z})^3 \times \mathbb{Z}/8 \mathbb{Z} $ & $192$ & $18$\\
$(\mathbb{Z}/4 \mathbb{Z} \times \mathbb{Z}/2 \mathbb{Z}) \rtimes \mathbb{Z}/2 \mathbb{Z}$ & $(\mathbb{Z}/2 \mathbb{Z})^3 \times \mathbb{Z}/8 \mathbb{Z} $ & $48$ & $32$\\
$(\mathbb{Z}/2 \mathbb{Z})^4 $ & $(\mathbb{Z}/2 \mathbb{Z})^{14}$ & $20160$ & $23$\\
\hline
\end{tabular}
\captionof{table}{}\label{tab:groups16}
\end{center}
Then, we compute the invariants for the $328$ classes that remain. The groups and the representatives of cohomology classes under the action of the automorphism group are labeled by numbers, in the same organization as in the previous cases. We already note that for the group $G_{11}$ and the cocycles $\omega_{11}^{14}$ and $\omega_{11}^{32}$, we get the same table. This phenomenon of a group giving rise to two non equivalent pointed categories that are (as we shall see) categorically Morita equivalent does not appear for groups of lower order.  We give here the set of Frobenius-Schur indicators and T-matrix for the center of this Morita class, and discuss the details later. In this table, we denote the primite $8$-th root of unity $\kappa:=\exp\left(\frac{i \pi}{4}\right)$.
\begin{center}
\begin{tabular}{*{5}{|c}|r|} \hline
$ \nu_0$ &$ \nu_1$ &$ \nu_2$ &$ \nu_4$ & $\theta$ & $m$\\
\hline
$1$ & $0$ & $1$ & $1$ & $1$ & $7$\\
\hline
$1$ & $1$ & $1$ & $1$ & $1$ & $1$\\
\hline
$2$ & $0$ & $-1$ & $1$ & $-1$ & $2$\\
\hline
$2$ & $0$ & $-1$ & $1$ & $1$ & $2$\\
\hline
$2$ & $0$ & $-1$ & $1$ & $-\kappa^2$ & $2$\\
\hline
$2$ & $0$ & $-1$ & $1$ & $\kappa^2$ & $2$\\
\hline
$2$ & $0$ & $-1$ & $1$ & $-\kappa$ & $2$\\
\hline
$2$ & $0$ & $-1$ & $1$ & $-\kappa^3$ & $2$\\
\hline
$2$ & $0$ & $-1$ & $1$ & $\kappa^3$ & $2$\\
\hline
$2$ & $0$ & $-1$ & $1$ & $\kappa$ & $2$\\
\hline
$2$ & $0$ & $1$ & $1$ & $-1$ & $2$\\
\hline
$2$ & $0$ & $1$ & $1$ & $1$ & $2$\\
\hline
$2$ & $0$ & $1$ & $1$ & $-\kappa^2$ & $2$\\
\hline
$2$ & $0$ & $1$ & $1$ & $\kappa^2$ & $2$\\
\hline
$2$ & $0$ & $1$ & $2$ & $-1$ & $2$\\
\hline
$2$ & $0$ & $1$ & $2$ & $1$ & $4$\\
\hline
$4$ & $0$ & $-1$ & $2$ & $-\kappa^2$ & $2$\\
\hline
$4$ & $0$ & $-1$ & $2$ & $\kappa^2$ & $2$\\
\hline
$4$ & $0$ & $1$ & $2$ & $-1$ & $2$\\
\hline
$4$ & $0$ & $1$ & $2$ & $1$ & $2$\\
\hline
\end{tabular}
\captionof{table}{FSI for $\mathcal{C}(G_7,\omega_7^{13})$, $\mathcal{C}(G_9\omega_9^1)$,
$\mathcal{C}(G_{11},\omega_{11}^{14})$ and $\mathcal{C}(G_{11}\omega_{11}^{32})$
}\label{tab:FSIquadruple16}
\end{center}

There are $230$ different tables among the $328$ that we computed. In detail, $160$ tables occur only once, $52$ twice, $12$ three times, three appear four times, two appear five times and finally one appears six times.

First, we deal with the cases when a table appears only twice. We use Corollary \ref{corcentral} like before. A new difference between the previous cases and this one is that this method does not answer the question whether the two considered categories are Morita equivalent in three cases, namely:
\begin{eqnarray}
(G_3,\omega_3^{14})& \text{ and } &(G_4, \omega_4 ^4)\text{ ,}\\
(G_3,\omega_3^{16})& \text{ and } &(G_4, \omega_4 ^8)\text{ ,}\\
(G_3,\omega_3^{17})& \text{ and } &(G_4, \omega_4 ^6)\text{ .}
\end{eqnarray}
Indeed, we could not find an extension $1 \rightarrow \hat{A} \rightarrow G_{3/4} \rightarrow H \rightarrow 1$ such that one of the six cocycles considered could been seen as the inflation of a $3$-cocycle on $H$. So we cannot use \cref{corcentral} and need  a different argument to conclude that the different couples define categorically Morita equivalent categories. More precisely, we need to consider a general group-theoretical category $\C(G,\omega,H,\mu)$, decide if it is pointed, and identify $(G',\omega')$ such that $\C(G,\omega,H,\mu)\cong\C(G',\omega')$. It turns out that for our purposes  it is sufficient to restrict ourselves to group-theoretical categories such that the $3$-cocycle is trivial on the subgroup, that is group-theoretical cateories $\mathcal{C}(G,\omega,H,1)$ such that $\omega(g,h,k)=1$ for all elements $g,h,k$ of $H$. This spares us another step in our treatment of cohomology groups (see \ref{sec:GAPcoho} above), but we note that it seems merely accidental that we find enough group-theoretical categories of this form. We can then arrange that the cocycle $\omega$ is adapted. To decide whether $\C(G,\omega,H,1)$ is pointed we can simply count its simples, or calculate the dimension of the simples in the form of their $0$-th Frobenius-Schur indicators (which is a reasonable only sincewe are interested in listing indicators anyway). If $\C(G,\omega,H,1)=\C(G',\omega')$ is thus found to be pointed, it is still not immediately obvious how to extract $G'$ and $\omega'$ (or rather the cohomology class of the latter). Here, we were successful with the following approach: We can calculate the Frobenius-Schur indicators for the group-theoretical category $\C(G,\omega,H,1)$ which is categorically Morita equivalent to $\C(G,\omega)$. If $\C(G,\omega,H,1)$ is to be equivalent to $\C(G',\omega')$, then the Frobenius-Schur indicators for $\C(G',\omega')$ need to coincide with the Frobenius-Schur indicators for $\C(G,\omega,H,1)$. Moreover $(G',\omega)$ is then categorically Morita equivalent to $\C(G,\omega)$, thus the invariants of the centers of $\C(G,\omega)$ and $\C(G',\omega')$ necessarily coincide (which is how we determined the pairs $(G,\omega)$ and $(G',\omega')$ that we are interested in in the first place). As it turns out, these conditions are sufficient to draw the desired conclusions.

First of all, we compute the Frobenius-Schur indicators for the $6$ considered pointed categories; it happens that for the three cases, the two pointed categories associated have different sets of indicators. So, to conclude it suffices to find a subgroup $H$ of $G_3$ such that some group-theoretical categories $\mathcal{C}(G_3, \omega_3^{14},H,\mu)$, $\mathcal{C}(G_3, \omega_3^{16},H,\mu)$ and $\mathcal{C}(G_3, \omega_3^{17},H,\mu)$, are pointed and have, respectively, the same sets of Frobenius-Schur indicators as the pointed categories $\mathcal{C}(G_4,\omega_4^4)$, $\mathcal{C}(G_4,\omega_4^8)$ and $\mathcal{C}(G_4,\omega_4^6)$. In order to obtain group-theoretical categories that are pointed, it is sufficient to consider abelian normal subgroups of $G_3$. Also, we restrict ourselves to subgroups $H$ such that the respective restrictions of the cocycles $\omega_3^{14}$, $\omega_3^{16}$ and $\omega_3^{17}$ are trivial in $C^3(H,\mathbb{C}^\times)$. Take $\delta \in \{ 14,16,17\}$. For each value of $\delta$ we find such an $H$, such that the group-theoretical category $\mathcal{C}(G_3,\omega_3^\delta,H,1)$ is pointed; we see this by computing the invariants first, and then look of the values of the $0$-th Frobenius-Schur indicators, which gives the dimensions of the simple objects. Moreover, it has the same indicators as the pointed category $\mathcal{C}(G_4,\omega_4^\delta)$. So, we can conclude for those three couples that the associated pointed fusion categories are categorically Morita equivalent, and so the doubles are gauge equivalent. 
\begin{remark}
We have used the fact that Frobenius-Schur indicators distinguish certain pointed categories (up to equivalence) to help us conclude that they are categorically Morita equivalent. Note that, in general, indicators do not classify pointed categories: The smallest example is for pointed categories of rank $8$, for which there exist $47$ non-equivalent categories, but only $34$ different sets of indicators.
\end{remark}
For the remaining $18$ tables that occur for more than two categories, except one, we can use the same kind of argument as for Table \ref{tab:FSItriple8} in \ref{sssec:8}. More precisely, for a given table, suppose that we find a group of the form $H \times Q$ that has only one cocycle affording this table (let us denote it $\mu$). Then for every other couple $(G,\omega)$ sharing this table, if we find a central extension $1 \rightarrow H \rightarrow G \rightarrow Q \rightarrow 1$ with $\omega$ the inflation of a $3$-cocycle on $Q$, the conditions of \cref{theo21} are satisfied and there exists a Morita equivalence between $\mathcal{C}(G,\omega)$ and $\mathcal{C}(H\times Q,\mu)$. This reasoning gives us only a part of the categorical Morita equivalence relation, but in our cases a simple inspection shows that the fact that categorical Morita equivalence is an equivalence relation implies that all couples sharing an identical table are in fact equivalent.

We illustrate this with an exemple. The six couples $(G_4, \omega_4^1)$, $(G_{10}, \omega_{10}^{6})$, $(G_{11}, \omega_{11}^{13})$,  $(G_{12}, \omega_{12}^{2})$,  $(G_{13}, \omega_{13}^{2})$ and  $(G_{14}, \omega_{14}^{8})$ share the same table. We will show by the above arguments that each of the pointed categories associated to the first five couples is categorically Morita equivalent to $\mathcal C(G_{14},\omega_{14}^8)$; note that $\omega_{14}^8$ is the only cocycle orbit on $G_{14}$ giving rise to the invariants in question on the double. The group $G_4$ is a semidirect product $\mathbb{Z}/4 \mathbb{Z} \rtimes \mathbb{Z}/4\mathbb{Z}$ and its center $H$ as well as the quotient $Q=G_4/H$ is isomorphic to $\mathbb{Z}/2 \mathbb{Z} \times \mathbb{Z}/2\mathbb{Z}$. The trivial cocycle on this quotient inflates to the trivial cocycle on $G_4$ and so we conclude an equivalence between the couples $(G_4, \omega_4^1)$ and $(G_{14}, \omega_{14}^{8})$, as $G_{14}$ is isomorphic to $H\times Q$.  For all of the four groups $G_{10}$, $G_{11}$, $G_{12}$ and $G_{13}$, there exists a central extension with subgroup isomorphic to $\mathbb{Z}/2 \mathbb{Z}$ and quotient isomorphic to $(\mathbb{Z}/2\mathbb{Z})^3$: Each of the four respective cocycles can be seen as inflations from cocycles on the respective quotients. Since also $G_{14}=\mathbb{Z}/2 \mathbb{Z}\times (\mathbb{Z}/2\mathbb{Z})^3$, we can conclude that the couple $(G_{14}, \omega_{14}^8)$ is equivalent to any of the five others, and so those six couples define Morita equivalent pointed categories. The procedure we have outlined in this particular example of six categorically Morita equivalent categories can be easily implemented in order to let the computer automatically search for such equivalences.

There is one single table for which the preceding argument does not say if the data that give this table determine equivalent categories. This table is the one we emphasized before and occurs for the couples:
\begin{equation}\label{quad}
(G_7,\omega_7^{13}) \text{ , } (G_9,\omega_9^1) \text{ , } (G_{11}, \omega_{11} ^{14})\text{ and } (G_{11}, \omega_{11} ^{32}) \text{ .}
\end{equation}
Again, we need to compute Frobenius-Schur indicators of the $4$ pointed categories and some group-theoretical categories. The two pointed categories $\mathcal{C}(G_{11}, \omega_{11} ^{14})$ and $\mathcal{C}(G_{11}, \omega_{11} ^{32})$ have the same set of Frobenius-Schur indicators, but it is different from the tables for $\mathcal{C}(G_{7}, \omega_{7} ^{13})$ and $\mathcal{C}(G_{9}, \omega_{9} ^{1})$; also, the two last ones are different from each other. We denote those $3$ sets of invariants respectively $PT_1$, $PT_2$ and $PT_3$. Then, we search for abelian normal subgroups of the groups $G_7$, $G_9$ and $G_{11}$, such that the restrictions of the involved cocycles are trivial.
\begin{center}
\begin{tabular}{*{8}{|c}||r|} \hline
$ \nu_0$ &$ \nu_1$ &$ \nu_2$ &$ \nu_3$ &$ \nu_4$ &$ \nu_5$ &$ \nu_6$ &$ \nu_7$ &$m$\\
\hline
\hline
$1$ & $0$ & $0$ & $0$ & $-1$ & $0$& $0$& $0$& $2$\\
\hline
$1$ & $0$ & $-1$ & $0$ & $1$ & $0$& $-1$& $0$& $5$\\
\hline
$1$ & $0$ & $0$ & $0$ & $1$ & $0$& $0$& $0$& $2$\\
\hline
$1$ & $0$ & $1$ & $0$ & $1$ & $0$& $1$& $0$& $6$\\
\hline
$1$ & $1$ & $1$ & $1$ & $1$ & $1$& $1$& $1$& $1$\\
\hline
\end{tabular}
\captionof{table}{PT1: FSI for $\mathcal{C}(G_{11}, \omega_{11} ^{14})$ and $\mathcal{C}(G_{11}, \omega_{11} ^{32})$}
\end{center}
\begin{center}
\begin{tabular}{*{8}{|c}||r|} \hline
$ \nu_0$ &$ \nu_1$ &$ \nu_2$ &$ \nu_3$ &$ \nu_4$ &$ \nu_5$ &$ \nu_6$ &$ \nu_7$ &$m$\\
\hline
\hline
$1$ & $0$ & $0$ & $0$ & $0$ & $0$& $0$& $0$& $4$\\
\hline
$1$ & $0$ & $-1$ & $0$ & $1$ & $0$& $-1$& $0$& $4$\\
\hline
$1$ & $0$ & $0$ & $0$ & $1$ & $0$& $0$& $0$& $2$\\
\hline
$1$ & $0$ & $1$ & $0$ & $1$ & $0$& $1$& $0$& $5$\\
\hline
$1$ & $1$ & $1$ & $1$ & $1$ & $1$& $1$& $1$& $1$\\
\hline
\end{tabular}
\captionof{table}{PT2: FSI for $\mathcal{C}(G_{7}, \omega_{7} ^{13})$}
\end{center}
\begin{center}
\begin{tabular}{*{8}{|c}||r|} \hline
$ \nu_0$ &$ \nu_1$ &$ \nu_2$ &$ \nu_3$ &$ \nu_4$ &$ \nu_5$ &$ \nu_6$ &$ \nu_7$ &$m$\\
\hline
\hline
$1$ & $0$ & $0$ & $0$ & $0$ & $0$& $0$& $0$& $4$\\
\hline
$1$ & $0$ & $0$ & $0$ & $1$ & $0$& $0$& $0$& $10$\\
\hline
$1$ & $0$ & $1$ & $0$ & $1$ & $0$& $1$& $0$& $1$\\
\hline
$1$ & $1$ & $1$ & $1$ & $1$ & $1$& $1$& $1$& $1$\\
\hline
\end{tabular}
\captionof{table}{PT3: FSI for $\mathcal{C}(G_{9}, \omega_{9} ^{1})$}
\end{center}
There is only one subgroup $H_1$ of $G_7$ such that $\omega_7^{13}|_{H_1 \times H_1 \times H_1}=1$. The category $\mathcal{C}(G_7,\omega_7^{13},H_1,1)$ is pointed, and its set of Frobenius-Schur indicators is $PT_3$. So, we can say that the category $\mathcal{C}(G_7,\omega_7^{13})$ is categorically Morita equivalent to at least one of the two categories $\mathcal{C}(G_{11} ,\omega_{11}^{14})$ and $\mathcal{C}(G_{11} ,\omega_{11}^{14})$, but we cannot directly know to which one of them.

We get much more information with the couple $(G_9,\omega_9^1)$: $\omega_9^1$ is the trivial cocycle on $G_9$, so each abelian normal subgroups of $G_9$ define a group-theoretical category; there are $3$ of them. Let us denote $H_2$ the normal abelian group of the quasi-dihedral group $G_9 \cong QD_{16}$ isomorphic to a cyclic group of order $8$. The category $\mathcal{C}(G_9,\omega_9^1,H_2,1)$ is pointed and its set of invariants is $PT_1$. So, we can conclude that the categories $\mathcal{C}(G_7,\omega_7^{13})$ and $\mathcal{C}(G_9,\omega_9^1)$ are categorically Morita equivalent. The two other subgroups give pointed categories that are equivalent to either $\mathcal{C}(G_{11} ,\omega_{11}^{14})$ or $\mathcal{C}(G_{11} ,\omega_{11}^{14})$, but again we could say no more.

Now, we have a look on the abelian normal subgroups of $G_{11}$: there are $13$ such groups, $\omega_{11}^{14}$ and $\omega_{11}^{32}$ are trivial on respectively $5$ and $3$ of them. We found subgroups $H_3$ and $H_4$, such the restrictions $\omega_{11}^{14} |_{H_3 \times H_3 \times H_3}$ and $\omega_{11}^{32} |_{H_4 \times H_4 \times H_4}$ are trivial and the categories\\ $\mathcal{C}(G_{11},\omega_{11}^{14},H_3,1)$ and $\mathcal{C}(G_{11},\omega_{11}^{32},H_4,1)$ are pointed and both have their set of indicators equal to $PT_1$. So we can conclude that the four couples in \cref{quad} define categorically Morita equivalent pointed fusion categories.

Finally, we are able to state:
\begin{theorem}\label{theo16}
There are exactly $230$ Morita equivalence classes of poinbted fusion categories of rank $16$.
\end{theorem}
\subsection{Results and perspectives}\label{ssec:resandpers}

In this section we gather the results obtained in the paper and ask a few questions that we think deserve to be observed. The results are either results of complete classification of pointed categories up to Morita equivalence or properties of completeness of the invariants used. Table \ref{tab:global} summarizes the number of Morita equivalence classes and usual equivalence classes of pointed categories for each order of $G$, from $2$ to $31$. It includes therefore \cref{theo8,theo27,theo16}, along with other particular orders we did not discuss in details in the present paper. Among them, groups of order $p$ and $p^2$ where $p$ is a prime show a rather predictable behavior; we will try to treat this in general in a separate paper. The table also contains the groups of order $p^3$ for $p=2,3$, as seen in the previous sections; here also, we think some common behavior could be extracted for all $p$, although this case should be more complicated. Also, some groups of order $pq$ where $p$ and $q$ are distinguished prime numbers are treated, but we will again look more closely at them in another study. For the proofs of \cref{completeclass} and \cref{completeinv} we only note that they are obtained by the same arguments and techniques used in section \ref{sec:work}. Except for the case $|G|$ discussed in \ref{sssec:16}, indicators of group-theoretical categories are necessary only for $|G|=18$.
\begin{theorem}\label{completeclass}
There are exactly $1126$ non-Morita equivalent pointed categories $\C(G,\omega)$ with $|G[ <32$.
\end{theorem}
\begin{theorem}\label{completeinv}
The set of Frobenius-Schur indicators and T-matrix, and \emph{a fortiori} the set of S and T-matrices, is a complete set of invariants for the modular categories of the form $\mathcal{Z}(\C(G,\omega))$ with $|G[ <32$.
\end{theorem}

\begin{center}
\begin{tabular}{|c|c|c|c|}
\hline
\rotatebox{90}{$|G|$} & \rotatebox{90}{Number of groups} & \rotatebox{90}{\begin{tabular}{c}Number of\\ non-equivalent\\ pointed categories\\ $\mathcal{C}(G,\omega)$ \end{tabular}} & \rotatebox{90}{\begin{tabular}{c}Number of\\ non-gauge\\ equivalent doubles\\ $D_\omega(G)$\end{tabular}}\\
\hline\hline
2 & 1 & 2 & 2\\ \hline
3 & 1 & 3 & 3\\ \hline
4 & 2 & 8 & 7\\ \hline
5 & 1 & 3 & 3\\ \hline
6 & 2 & 12 & 12\\ \hline
7 & 1 & 3 & 3\\ \hline
8 & 5 & 47 & 38\\ \hline
9 & 2 & 10 & 9\\ \hline
10 & 2 & 12 & 12\\ \hline
11 & 1 & 3 & 3\\ \hline
12 & 5 & 60 & 54\\ \hline
13 & 1 & 3 & 3\\ \hline
14 & 2 & 12 & 12\\ \hline
15 & 2 & 9 & 9\\ \hline
16 & 14 & 328 & 230\\ \hline
17 & 1 & 3 & 3\\ \hline
18 & 5 & 58 & 54\\ \hline
19 & 1 & 3 & 3\\ \hline
20 & 5 & 58 & 52\\ \hline
21 & 2 & 12 & 12\\ \hline
22 & 2 & 12 & 12\\ \hline
23 & 1 & 3 & 3\\ \hline
24 & 15 & 474 & 393 \\ \hline
25 & 2 & 10 & 9\\ \hline
26 & 2 & 12 & 12\\ \hline
27 & 5 & 61 & 47\\ \hline
28 & 4 & 54 & 48\\ \hline
29 & 1 & 3 & 3\\ \hline
30 & 4 & 72 & 72\\ \hline
31 & 1 & 3 & 3\\ \hline
\end{tabular}
\captionof{table}{Number of pointed categories and twisted quantum doubles for 1<$|G|$<32}\label{tab:global}
\end{center}
Not surprisingly, we have tried the methods we have described above to treat pointed categories $\mathcal{Z}(\C(G,\omega))$ with $|G|=32$. We postpone details about this class of pointed fusion categories, as we have not yet obtained the exact number of non-equivalent Morita classes. We just note that there are at most $4081$ such classes, this number being the sum of automorphism orbits of cocycles running over the finite groups of order $32$, and at least $2315$, corresponding to the number of different sets of Frobenius-Schur indicators and ribbon twist of the centers of these categories. In order to get over the difficulties we encounter in this case, we think that explicit computation of Frobenius-Schur indicators of general group-theoretical categories cannot be avoided. However, we also find $20$ sets of Frobenius-Schur indicators and ribbon twists such that each occurs for exactly two couples of group-cocycle pairs $(G_1,\omega_1)$ and $(G_2,\omega_2)$ with $G_1=G_2$ and such that the corresponding pointed categories  $\C(G_1,\omega_1)$ and $\C(G_2,\omega_2)$ are not distinghuished by their own Frobenius-Schur indicators. Computing indicators of group-theoretical categories is then useless for these cases, and we are not able, using only the methods above, to answer if the categories $\mathcal{Z}(\C(G_1,\omega_1))$ and $\mathcal{Z}(\C(G_2,\omega_2))$ are indeed equivalent.

In a separate paper we will present examples of two modular categories of the form $\mathcal{Z}(\C(G,\omega))$ which are non-equivalent but share the same the set of Frobenius-Schur indicators and ribbon twist. Moreover, we prove that even the set of S and T-matrices is not a complete invariant in this case, providing a counterexample for the ``belief'' professed in \cite{Bruillard}. The examples we found are for the semi-direct products $\ZZ/q\ZZ \semdir \ZZ/p\ZZ$ where $p$ and $q$ are two primes greater than $3$ such that $q|p-1$; the smallest such example thus occurs for the nonabelian group of order $55$, while we have seen that up to order $31$ categorical Morita equivlence classes of pointed categories are distinguished by the Frobenius-Schur indicators and T-matrices of their centers. In particular, up to dimension $31$ the centers of pointed fusion categories are distinguished by their modular data, but not in dimension $55$. It may be that the smallest example occurs for some order smaller than $55$, but we were not able to decide this.

\bibliographystyle{habbrv}

\bibliography{biblio.bib}
\end{document}